\documentclass[preprint,sort&compress,number]{elsarticle}

\usepackage{amsmath,amsthm,amssymb,latexsym,amscd}
\usepackage{geometry}
\usepackage{algorithm}
\usepackage{algpseudocode}
\usepackage{dsfont}
\usepackage{color}
\usepackage{lineno}
\usepackage{changes}
\usepackage[titletoc]{appendix}
\usepackage{url}
\usepackage{amsmath}
\usepackage{bm}
\usepackage{textgreek}
\usepackage{accents}
\usepackage{multirow}
\usepackage{subfigure}
\usepackage{caption}
\usepackage[bb=boondox]{mathalfa}
\usepackage{hyperref}

\newcommand{\ATone}{\texttt{AT$_1$}}
\newcommand{\ATtwo}{\texttt{AT$_2$}}

\newcommand{\mrm}{\mathrm}

\newcommand{\mathd}{\mathrm{d}}

\newcommand{\RNum}[1]{\uppercase\expandafter{\romannumeral #1\relax}}

\DeclareMathOperator{\tr}{tr}
\DeclareMathOperator{\dev}{dev}

\algtext*{EndIf}
\algtext*{EndWhile}
\algtext*{EndFor}

\newtheorem{remark}{Remark}

\definechangesauthor[name=VZR,color=brown]{VZR}

\makeatletter
\def\ps@pprintTitle{%
	\let\@oddhead\@empty
	\let\@evenhead\@empty
	\def\@oddfoot{}%
	\let\@evenfoot\@oddfoot}
\makeatother

\title{Orthogonal decomposition of anisotropic constitutive models for the phase field approach to fracture}
\author[UFZ]{Vahid Ziaei-Rad}
\author[UFZ]{Mostafa Mollaali}
\author[TUBAF]{Thomas Nagel}
\author[UFZ,TUD]{Olaf Kolditz}
\author[UFZ]{Keita Yoshioka}
\address[UFZ]{Helmholtz Centre for Environmental Research -- UFZ, Leipzig, Germany}
\address[TUBAF]{Geotechnical Institute, Technische Universität Bergakademie Freiberg, Germany}
\address[TUD]{Technische Universität Dresden, Germany}

\begin{document}
    \begin{abstract}
        We propose a decomposition of constitutive relations into crack-driving and persistent portions, specifically designed for materials with anisotropic/orthotropic behavior in the phase field approach to fracture to account for the tension-compression asymmetry. This decomposition follows a variational framework, satisfying the orthogonality condition for anisotropic materials. This implies that the present model can be applied to arbitrary anisotropic elastic behavior in a three-dimensional setting. On this basis, we generalize two existing models for tension-compression asymmetry in isotropic materials, namely the `volumetric-deviatoric' model \cite{Amor2009} and the `no-tension' model \cite{Freddi2010}, towards materials with anisotropic nature. 
        Two benchmark problems, single notched tensile  shear tests, are used to study the performance of the present model. The results can retain the anisotropic constitutive behavior and the tension-compression asymmetry in the crack response, and are qualitatively in accordance with the expected behavior for orthotropic materials. Furthermore, to study the direction of maximum energy dissipation, we modify the surface integral based energy release computation, $G_\theta$, to account only for the crack-driving energy. The computed energies with our proposed modifications predict the fracture propagation direction correctly compared with the standard $G_\theta$ method.
    \end{abstract}
    \maketitle
	\section{Introduction}


    The phase field method has steadily gained popularity for over a decade because of its convenience in simulating complex fracture processes, including crack nucleation, propagation, branching and merging~\cite{Bourdin2019}. The model was originally conceived as the variational formulation of brittle fracture by Francfort and Marigo~\cite{Francfort1998}, and its regularization by Bourdin et al~\cite{Bourdin2000}. The formulation solves crack problems by minimizing an energy functional that consists of the elastic energy and the crack surface energy. This significantly reduces the implementation complexities as the crack evolution is a natural outcome of the solution. Therefore, there is no need for a crack tracking algorithm or additional criteria for crack branching and merging. The method has found its extension in subjects such as composite failure~\cite{baldelli2014variational,Alessi2017, Yoshioka2021int}, ductility~\cite{alessi2014gradient, Ambati2015,Kuhn2016, miehe2016ductile,you2021brittle}, dynamic fracture~\cite{Bourdin2011, Borden2012,Nguyen2018, Yin2020}, hydraulic fracture~\cite{Bourdin2012, Wheeler2014, Wilson2016, Heider2017,zhou2018phase}, and environment assisted failure~\cite{martinez2018phase,schuler2020chemo} to name a few. \newline
    When mechanics of natural materials is concerned, especially geomaterials, they often exhibit anisotropic behaviors~\cite{carcione1995constitutive, Amadei1996,parisio2018formulation}.
    Anisotropic material response manifests itself in terms of deformation~\cite{nixon2010anisotropic, heng2015experimental} and/or preferential fracture propagation~\cite{wang2019anisotropic,li2021research}.
    In modeling anisotropy within the phase field framework, anisotropic fracture toughness has been mostly achieved through manipulation of the fracture surface energy term where one applies structural tensors~\cite{li2015phase,teichtmeister2017phase,bleyer2018phase,Li2019aniso,bilgen2021phase,gerasimov2022second} or assign direction dependent fracture toughness~\cite{bryant2018mixed,noii2020adaptive,fei2021double,rezaei2022anisotropic,gerasimov2022second,ulloa2022micromechanics} to encourage fracture propagation in certain directions. 
    While many studies are available for anisotropic fracture toughness, only a few have studied anisotropic constitutive model in phase field approaches~\cite{Nguyen2020, DIJK2020, luo2022phase}. \newline
    Another material response to consider is tension-compression asymmetry in failure, i.e., a crack is less likely to propagate under compressive loadings. Consequently, there has to be a unilateral constraint in material degradation. Otherwise, unphysical crack propagation is predicted under compressive loading. How to account for this unilateral constraint is one of the main challenges in phase field modeling. 
    Note that such a model is sometimes referred to in the literature as an `anisotropic' phase field model \cite{Ambati2015}, however, in this work we use the term of `tension-compression asymmetry' for this effect (see also \cite{Moerman2016}), and preserve `anisotropy' to only address a material whose constitutive behavior is, unlike isotropic materials, directionally dependent. \newline
    With the need of tension-compression asymmetry in material degradation, several existing approaches consist of additively partitioning the sound elastic energy density $\varPsi$ into two portions, a crack-driving portion $\varPsi^\text{t}$ and a persistent portion $\varPsi^\text{c}$~\cite{Amor2009,Miehe2010_variational, Freddi2010,Steinke2018,wu2018length,de2021nucleation}.
    If we are able to partition the strain energy density based on the additive decomposition of the strain tensor, that is, 
    \begin{equation*}
        \begin{aligned}
            \bm{\varepsilon} = \bm{\varepsilon}^\text{t} + \bm{\varepsilon}^\text{c}, \quad \varPsi(\bm{\varepsilon}) = \varPsi^\text{t}(\bm{\varepsilon}^\text{t}) + \varPsi^\text{c}(\bm{\varepsilon}^\text{c}), \\
            \varPsi^\text{t}(\bm{\varepsilon}^\text{t}) = \frac12 \mathbb{C} \bm{\varepsilon}^\text{t}\cdot\bm{\varepsilon}^\text{t}, \quad 
            \varPsi^\text{c}(\bm{\varepsilon}^\text{t}) = \frac12 \mathbb{C} \bm{\varepsilon}^\text{c}\cdot\bm{\varepsilon}^\text{c}, 
        \end{aligned}
    \end{equation*}
    the constitutive behaviors are characterized by the same elasticity tensor $\mathbb{C}$ where a dot $\cdot$ denotes the inner products between the two vectors or tensors of the same order. Then, there exists a local \emph{variational} principle from which several existing models can be derived:
    \begin{equation}\label{Eq:variational_principle}
        \begin{aligned}
            \bm{\varepsilon}^\text{t} = \underset{\bm{e}^\text{t}\in\mathcal{S}^\text{t}}{\mrm{argmin}} \, 
            \mathbb{C}\bigl(\bm{\varepsilon} - \bm{e}^\text{t}\bigr)\cdot\bigl(\bm{\varepsilon} - \bm{e}^\text{t}\bigr),
        \end{aligned}
    \end{equation}
    where the strain that contributes to damage ($\bm{\varepsilon}^\text{t}$) is defined as the orthogonal projection of the strain tensor $\bm{\varepsilon}$ onto the convex space $\mathcal{S}^\text{t}$ with respect to the energy norm defined by $\mathbb{C}$. Therefore, the modeling of tension-compression asymmetry in fracture response is reduced to the identification of such convex cone $\mathcal{S}^\text{t}$ to represent $\bm{\varepsilon}^\text{t}$~\cite{Freddi2010}. Also, $\bm{\varepsilon}^\text{c}$ is in a polar cone $\mathcal{S}^\text{c}=\Big\{\bm{e}^\text{c}\Big|\bm{e}^\text{c}\cdot\bm{e}^\text{t}\leq 0,~\forall\bm{e}^\text{t} \in \mathcal{S}^\text{t}\Big\}$. This way, we keep the variational structure of the crack problem by following the energy minimization. This is important because otherwise, one needs to introduce {ad hoc} criteria for the phase field problem. For more detailed analysis, we refer the readers to \cite{Marigo2016} and \cite{Freddi2010}. \newline
    \emph{A challenges in anisotropic materials} is to satisfy this variational structure.
    From convex analysis, $\bm{\varepsilon}^\text{t}$ can be equivalently characterized by
    \begin{equation*}
        \begin{aligned}
            -\mathbb{C}(\bm{\varepsilon} - \bm{\varepsilon}^\text{t})\cdot\bigl(\bm{e}^\text{t} - \bm{\varepsilon}^\text{t}\bigr) \geq 0,
            \quad \forall \bm{e}^\text{t}\in\mathcal{S}^\text{t}.
        \end{aligned}
    \end{equation*}
    Testing this with $\bm{e}^\text{t}=2\bm{\varepsilon}^\text{t}$ and $\bm{e}^\text{t}=\frac{1}{2} \bm{\varepsilon}^\text{t}$ yields the following orthogonality condition~\cite{Li2016grad}:
    \begin{equation}\label{Eq:orthogonality}
        \begin{aligned}
            \mathbb{C}\bm{\varepsilon}^\text{t}\cdot\bm{\varepsilon}^\text{c} = \mathbb{C}\bm{\varepsilon}^\text{c}\cdot\bm{\varepsilon}^\text{t} = 0,
        \end{aligned}
    \end{equation}
    which is the necessary and sufficient condition for the strain energy density to be decomposed based on the partition of the strain tensor. However, for anisotropic linear elastic materials the principal directions of the stress and strain no longer coincide, i.e., the orthogonality condition is not generally satisfied:
    \begin{equation*}
        \begin{aligned}
            \mathbb{C}\bm{\varepsilon}^\text{t}\cdot\bm{\varepsilon}^\text{c} \neq 0 \text{~or~} \mathbb{C}\bm{\varepsilon}^\text{c}\cdot\bm{\varepsilon}^\text{t} \neq 0.
        \end{aligned}
    \end{equation*}
    {Therefore, further treatments are required in strain decomposition that keep the orthogonality condition for anisotropic materials.} \newline
    Recently, He and Shao \cite{He-2019} proposed an elastic energy preserving transformation so that an orthogonal decomposition of the strain tensor can be achieved in the transformed space. Therein, a rigorously orthogonal spectral decomposition of the strain tensor is performed.
    Thus, with the orthogonality condition being satisfied in the new space, the strain energy density of an anisotropic material is respectively partitioned into a crack-driving portion and a persistent portion. Moreover, to decompose the strain tensor in a two-dimensional setting, a set of closed-form coordinate-free expressions are provided.
    Nguyen et al. \cite{Nguyen2020} directly employed these expressions of strain decomposition into the phase field formulation, and compared the results with volumetric-deviatoric and spectral-decomposition for isotropic materials.
    Van Dijk et al. \cite{DIJK2020} proposed a generalization of {volumetric-deviatoric} and spectral-decomposition to account for anisotropy of materials. Therein, the strain energy density is described as an isotropic tensor function of a type of tensorial square roots of the strain energy density.
    \newline
    In this work, we present generalized formulations of two existing models for tension-compression asymmetry, a model by Amor et al. \cite{Amor2009} (henceforth referred to as {volumetric-deviatoric}) and a model by Freddi and Royer \cite{Freddi2010} (henceforth referred to as {no-tension}), to account for materials with anisotropic constitutive behavior in a three-dimensional setting. We borrow the basic concepts from \cite{He-2019} and perform the energy preserving transformation so that the orthogonality condition is preserved for a decomposition of an arbitrary anisotropic constitutive model. Then, we present the implementation details of {volumetric-deviatoric} and {no-tension} accounting for the anisotropy in Section \ref{Sec:problem-statement}. In addition, we provide a set of numerical examples to illustrate the capabilities of the present method in Section \ref{Sec:numerical-examples}. The results provide a step forward when developing phase field fracture theories for brittle materials with an anisotropic nature and highlight the importance of a proper decomposition of the strain energy density. \newline
    \section{Problem statement}\label{Sec:problem-statement}
    This section introduces the theory and implementation details to account for the tension-compression asymmetry in the crack response of anisotropic materials. In \ref{Sec:problem-PF} we recapitulate the basic ingredients of a phase field formulation. Then, in \ref{Sec:problem-transformation} we perform an elastic energy preserving transformation so that an orthogonal decomposition of the strain tensor, and its consequent partitioning of the strain energy density, can be realized in the transformed space. Thereafter, in \ref{Sec:problem-unilateral} we pick the two existing models for the tension-compression asymmetry, {volumetric-deviatoric} and {no-tension}, and adapt them for anisotropic materials. Finally, in \ref{Sec:problem-anisotropy} we provide various types of the elasticity tensors that can be adapted into the model. Readers interested in directly using the model can jump to \eqref{Eq:D_tensors}, \eqref{Eq:strain_decomp_voldev}, \eqref{Eq:sigma}, \eqref{Eq:C_decomp} for anisotropic {volumetric-deviatoric} and \eqref{Eq:proj-dir-dev}, \eqref{Eq:strain_decomp_masonry}, \eqref{Eq:sigma_masonry_tensor}, \eqref{Eq:C_masonry_tensor} for anisotropic {no-tension}.
\subsection{Phase field formulation}\label{Sec:problem-PF}
The phase field formulation for brittle fracture is essentially based on minimization of the following functional
\begin{equation}\label{Eq:pi}
    \begin{aligned}
        \Pi_{\ell}(\boldsymbol{u},d) := \int_{\Omega} \varPsi\bigl[\boldsymbol{\varepsilon}(\boldsymbol{u}),d\bigr]\,\mathrm{d} \Omega - \int_{\partial\Omega_t } \boldsymbol{t}_{N} \cdot \boldsymbol{u} \,\mathrm{d}A - \int_{\Omega} \boldsymbol{b} \cdot \boldsymbol{u}\,\mathrm{d} \Omega +  \int_{\Omega} G_c \gamma (d, \nabla d) ~\mathrm{d}\Omega,
    \end{aligned}
\end{equation}
where $d:\Omega\rightarrow[0,1]$ represents the phase field, characterizing the material at its pristine state with $d = 0$ and at the fully damaged state with $d = 1$. Also, $\boldsymbol{t}_{N}:\partial\Omega_t\rightarrow\mathbb{R}^n$ is the prescribed traction boundary condition{,} and $ \boldsymbol{b}:\Omega\rightarrow\mathbb{R}^n$ is the body force. The strain energy density takes the following form 
\begin{equation}\label{Eq:strain_energy_density}
    \begin{aligned}
        \varPsi\bigl[\bm{\varepsilon},d\bigr] = g(d)\varPsi^\text{t}(\bm{\varepsilon}^\text{t}) + \varPsi^\text{c}(\bm{\varepsilon}^\text{c}),
    \end{aligned}
\end{equation}
where $g(d)=(1-d)^2$ is a degradation function. 
Furthermore, $ G_c > 0 $ is the critical crack energy release rate, and $\gamma(d,\nabla d)$ is the crack surface density per unit volume,
\begin{equation*}
     \gamma(d,\nabla d) = 
     \frac{1}{4c_w}\left( \frac{w(d)}{\ell}+\ell |\nabla d|^2 \right),
\end{equation*}
where $ \ell > 0 $ is the regularization length scale parameter which controls the width of the transition region of the smoothed crack. Crack geometric function $w(d)$ and normalization constant $c_w = \int_0^1 \sqrt{w(d)}\,\mathrm{d}d $ are model dependent. Specifically, for brittle fracture, classical examples are $w(d)=d^2$ and $c_w=1/2$ for the \ATtwo{} model; and $w(d)=d$ and $c_w=2/3$ for the \ATone{} model \cite{Pham2011Gradient}. \newline
In the sequel we derive the first and second variations of $\Pi_{\ell}$ as in \eqref{Eq:pi}, which will be needed for the discretized formulation. 
\paragraph{First variation of $\Pi_{\ell}$} Taking the first variation yields:
\begin{equation}\label{Eq:first_variation}
    \begin{aligned}
        \delta\Pi[(\bm{u}, d); (\Bar{\bm{u}}, \Bar{d})] &:= \int_{\Omega} 
        \bm{\sigma}\bigl[\bm{\varepsilon}(\bm{u}), d\bigr]
        \cdot\bm{\varepsilon}(\Bar{\bm{u}})\,\mathrm{d} \Omega 
        - \int_{\partial\Omega_t } \boldsymbol{t}_{N} \cdot \Bar{\boldsymbol{u}} \,\mathrm{d}A 
        - \int_{\Omega} \boldsymbol{b} \cdot \Bar{\boldsymbol{u}}\,\mathrm{d} \Omega \\
        &+ \int_{\Omega} g^{\prime}(d) \varPsi^\text{t}(\boldsymbol{\varepsilon})\Bar{d} ~\mathrm{d}\Omega 
        + \frac{G_c}{4c_w} \int_{\Omega} \Bigl(\frac{w^{\prime}(d)\Bar{d}}{\ell} + 2\ell \nabla d \cdot \nabla \Bar{d} \Bigr) ~\mathrm{d}\Omega,
    \end{aligned}
\end{equation}
where 
\begin{equation}\label{Eq:sigma_degraded}
    \begin{aligned}
        \bm{\sigma} := \frac{\partial\varPsi}{\partial\bm{\varepsilon}} = g(d)\bm{\sigma}^\text{t} + \bm{\sigma}^\text{c}, \quad \bm{\sigma}^\text{t} = \frac{\partial\varPsi^\text{t}}{\partial\bm{\varepsilon}}, \quad \bm{\sigma}^\text{c} = \frac{\partial\varPsi^\text{c}}{\partial\bm{\varepsilon}}.
    \end{aligned}
\end{equation}
\paragraph{Second variation of $\Pi_{\ell}$} We take another variation from \eqref{Eq:pi}:
\begin{equation}\label{Eq:second_variation}
    \begin{aligned}
        \delta^2\Pi[(\bm{u}, d); (\Bar{\bm{u}}, \Bar{d}); (\delta\bm{u}, \delta d)] &:= \int_{\Omega} 
        \bm{\varepsilon}(\delta\bm{u})\cdot\mathbb{C}\bigl[\bm{\varepsilon}(\bm{u}), d \bigr] \cdot \bm{\varepsilon}(\Bar{\bm{u}})\,\mathrm{d} \Omega 
        + \int_{\Omega} \bm{\varepsilon}(\delta\bm{u})\cdot \bm{\sigma}^\text{t}\bigl[\bm{\varepsilon}(\bm{u})\bigr] g^{\prime}(d)\Bar{d}\,\mathrm{d} \Omega \\
        &+ \int_{\Omega} \delta d g^{\prime}(d) \bm{\sigma}^\text{t}\bigl[\bm{\varepsilon}(\bm{u})\bigr]\cdot \bm{\varepsilon}(\Bar{\bm{u}}) \,\mathrm{d} \Omega
        + \int_{\Omega} \delta d g^{\prime\prime}(d)\varPsi^\text{t}(\bm{\varepsilon}) \Bar{d} \,\mathrm{d} \Omega \\
        &+ \frac{g_c}{4c_w} \int_{\Omega} \Biggl(\frac{\delta d w^{\prime\prime}(d) \Bar{d}}{\ell} + 2\ell\nabla(\delta d) \cdot \nabla\Bar{d} \Biggr) \,\mathrm{d} \Omega,
    \end{aligned}
\end{equation}
where 
\begin{equation}\label{Eq:C_degraded}
    \begin{aligned}
        \mathbb{C} := \frac{\partial^2\varPsi}{\partial\bm{\varepsilon}^2} = g(d) \mathbb{C}^\text{t} + \mathbb{C}^\text{c}, \quad \mathbb{C}^\text{t} = \frac{\partial^2\varPsi^\text{t}}{\partial\bm{\varepsilon}^2}, \quad \mathbb{C}^\text{c} = \frac{\partial^2\varPsi^\text{c}}{\partial\bm{\varepsilon}^2}.
    \end{aligned}
\end{equation}
The components of \eqref{Eq:strain_energy_density}, \eqref{Eq:sigma_degraded}, \eqref{Eq:C_degraded} are particularly obtained in the sequel.
\subsection{Transformation for anisotropic $\mathbb{C}$}\label{Sec:problem-transformation}
Let $\mathcal{S}$ and $\mathcal{S}^*$ be the set of all second-order strain and stress tensors, respectively. To enforce the orthogonality condition \eqref{Eq:orthogonality} for an arbitrary anisotropic $\mathbb{C}$, following He and Shao \cite{He-2019}, we define the transformed strain and stress spaces as
\begin{subequations}\label{Eq:all_transf}
    \begin{align}
        \tilde{\mathcal{S}}=\Big\{\tilde{\bm{\varepsilon}}\Big|\tilde{\bm{\varepsilon}}=\mathbb{C}^{1/2}\bm{\varepsilon},~\bm{\varepsilon}\in\mathcal{S}\Big\}, \label{Eq:eps_transf}\\
        {\tilde{\mathcal{S}}}^*=\Big\{\tilde{\bm{\sigma}}\Big|\tilde{\bm{\sigma}}=\mathbb{S}^{1/2}\bm{\sigma},~\bm{\sigma}\in\mathcal{S}^*\Big\}, \label{Eq:sigma_transf}
    \end{align}
\end{subequations}
so that
\begin{equation}\label{Eq:psi_tilde_decomp}
    \begin{aligned}
        \varPsi(\bm{\varepsilon}) = \frac12 \tilde{\bm{\varepsilon}}\cdot\tilde{\bm{\varepsilon}} = \frac12 \tilde{\bm{\sigma}}\cdot\tilde{\bm{\sigma}}.
    \end{aligned}
\end{equation}
{Note that $\mathbb{C}^{1/2}$ in \eqref{Eq:eps_transf} is the square root of $\mathbb{C}$, and $\mathbb{S}^{1/2}$ in \eqref{Eq:sigma_transf} is the square root of the compliance tensor $\mathbb{S}=\mathbb{C}^{-1}$, see Remark \eqref{Re:C_squar} for how to compute $\mathbb{C}^{1/2}$.}
Now let $\tilde{\mathcal{S}}$ and ${\tilde{\mathcal{S}}}^*$ be split into two convex subsets such that
\begin{subequations}\label{Eq:S-st-split}
    \begin{align}
        \tilde{\mathcal{S}}={\tilde{\mathcal{S}}}^{t} \cup {\tilde{\mathcal{S}}}^{c}, \label{Eq:S-st-split-strain} \\
        {\tilde{\mathcal{S}}}^*={\tilde{\mathcal{S}}}^{*t} \cup {\tilde{\mathcal{S}}}^{*c}. \label{Eq:S-st-split-stress}
    \end{align}
\end{subequations}
\eqref{Eq:S-st-split-strain} means that any element $\tilde{\bm{\varepsilon}}\in\tilde{\mathcal{S}}$ is decomposed into a crack-driving portion $\tilde{\bm{\varepsilon}}^{t}$ and a persistent portion $\tilde{\bm{\varepsilon}}^{c}$ such that
\begin{equation}\label{Eq:eps_tild_cond}
    \begin{aligned}
        \tilde{\bm{\varepsilon}} = \tilde{\bm{\varepsilon}}^\text{t} + \tilde{\bm{\varepsilon}}^\text{c}, \quad
        \tilde{\bm{\varepsilon}}^\text{t}\cdot\tilde{\bm{\varepsilon}}^\text{c}=0.
    \end{aligned}
\end{equation}
Similarly for the stress, one writes
\begin{equation*}
    \begin{aligned}
        \tilde{\bm{\sigma}} = \tilde{\bm{\sigma}}^\text{t} + \tilde{\bm{\sigma}}^\text{c}, \quad
        \tilde{\bm{\sigma}}^\text{t}\cdot\tilde{\bm{\sigma}}^\text{c}=0.
    \end{aligned}
\end{equation*}
Thus, with reference to \eqref{Eq:psi_tilde_decomp} the strain energy in the transformed spaces $\tilde{\mathcal{S}}$ or ${\tilde{\mathcal{S}}}^*$ can be decomposed as:
\begin{equation}\label{Eq:psi_tild_split}
    \begin{aligned}
        \varPsi(\tilde{\bm{\varepsilon}}) = \frac12\tilde{\bm{\varepsilon}}^\text{t}\cdot\tilde{\bm{\varepsilon}}^\text{t} + \frac12 \tilde{\bm{\varepsilon}}^\text{c}\cdot\tilde{\bm{\varepsilon}}^\text{c}, \quad
        \varPsi(\tilde{\bm{\sigma}}) = \frac12\tilde{\bm{\sigma}}^\text{t}\cdot\tilde{\bm{\sigma}}^\text{t} + \frac12 \tilde{\bm{\sigma}}^\text{c}\cdot\tilde{\bm{\sigma}}^\text{c}.
    \end{aligned}
\end{equation}
With $\tilde{\bm{\varepsilon}}^\text{t}$ and $\tilde{\bm{\varepsilon}}^\text{c}$ obtained, one can determine the inverse transformations:
\begin{equation}\label{Eq:tilde_ortho_cond}
    \begin{aligned}
        {\boldsymbol{\varepsilon}}^{t}=\mathbb{C}^{-1/2}\tilde{\boldsymbol{\varepsilon}}^{t}, \quad {\boldsymbol{\varepsilon}}^{c}=\mathbb{C}^{-1/2}\tilde{\boldsymbol{\varepsilon}}^{c}.
    \end{aligned}
\end{equation}
And similarly for the stress:
\begin{equation*}
    \begin{aligned}
        \bm{\sigma}^\text{t} = \mathbb{S}^{-1/2}\tilde{\bm{\sigma}}^\text{t}, \quad \bm{\sigma}^\text{c} = \mathbb{S}^{-1/2}\tilde{\bm{\sigma}}^\text{c}.
    \end{aligned}
\end{equation*}
One can observe that orthogonality also holds for the quantities transformed back to the physical stress and strain spaces.
Henceforth, the convex subsets are summarized as follows:
\begin{subequations}\label{Eq:S-tc-remap}
    \begin{align}
        \mathcal{S}^\text{t}=\Big\{\bm{\varepsilon}^\text{t}\Big|\bm{\varepsilon}^\text{t}=\mathbb{C}^{-1/2}\tilde{\bm{\varepsilon}}^\text{t},~\tilde{\bm{\varepsilon}}^\text{t}\in{\tilde{\mathcal{S}}}^{t}\Big\}, \quad \mathcal{S}^\text{c}=\Big\{\bm{\varepsilon}^\text{c}\Big|\bm{\varepsilon}^\text{c}=\mathbb{C}^{-1/2}\tilde{\bm{\varepsilon}}^\text{c},~\tilde{\bm{\varepsilon}}^\text{c}\in{\tilde{\mathcal{S}}}^{c}\Big\}, \label{Eq:S-tc-remap-strain} \\
        \mathcal{S}^{*t}=\Big\{\bm{\sigma}^\text{t}\Big|\bm{\sigma}^\text{t}=\mathbb{S}^{-1/2}\tilde{\bm{\sigma}}^\text{t},~\tilde{\bm{\sigma}}^\text{t}\in{\tilde{\mathcal{S}}}^{*t}\Big\}, \quad
        \mathcal{S}^{*c}=\Big\{\bm{\sigma}^\text{c}\Big|\bm{\sigma}^\text{c}=\mathbb{S}^{-1/2}\tilde{\bm{\sigma}}^\text{c},~\tilde{\bm{\sigma}}^\text{c}\in{\tilde{\mathcal{S}}}^{*c}\Big\}. \label{Eq:S-tc-remap-stress}
    \end{align}
\end{subequations}
\begin{remark}\label{Re:C_squar}
We derive the square root of $\mathbb{C}$ as follows. As $\mathbb{C}$ has major symmetry and is positive-definitive, one can write its spectral decomposition as below:
\begin{equation*}
    \begin{aligned}
        \mathbb{C}=\underset{i}{\sum}\Lambda_i~\textbf{\textomega}_i\otimes\textbf{\textomega}_i,
    \end{aligned}
\end{equation*}
where $\Lambda_i>0$ are the eigenvalues of $\mathbb{C}$ and \textbf{\textomega}$_i$ are the second-order orthonormal eigentensors associated to $\Lambda_i$, such that $\textbf{\textomega}_i\cdot\textbf{\textomega}_j=\delta_{ij}$.
Then, we can compute the square root of $\mathbb{C}$ by:
\begin{equation}\label{Eq:C_spectral_decomp}
    \begin{aligned}
        \mathbb{C}^{1/2}=\underset{i}{\sum}\Lambda_i^{1/2}\textbf{\textomega}_i\otimes\textbf{\textomega}_i.
    \end{aligned}
\end{equation}
The same procedure can be followed for the compliance tensor $\mathbb{S}$.
\end{remark}
\paragraph*{Summary of \ref{Sec:problem-transformation}} {We aimed at satisfying the orthogonality condition \eqref{Eq:orthogonality} for anisotropic $\mathbb{C}$. We illustrated that for any $\bm{\varepsilon}\in\mathcal{S}$, there exist $\bm{\varepsilon}^\text{t}\in\mathcal{S}^{t}$ and $\bm{\varepsilon}^\text{c}\in\mathcal{S}^{c}$ such that the orthogonality condition holds. This was achieved through an energy preserving transformation ($\varPsi(\bm{\varepsilon}) = \varPsi(\tilde{\bm{\varepsilon}})$) to a new space $\tilde{\bm{\varepsilon}}\in\tilde{\mathcal{S}}$ where $\varPsi(\tilde{\bm{\varepsilon}}) = \varPsi(\tilde{\bm{\varepsilon}}^\text{t}) + \varPsi(\tilde{\bm{\varepsilon}}^\text{c})$. The remaining step is to make a choice for $\tilde{\mathcal{S}}^\text{t}$ and $\tilde{\mathcal{S}}^\text{c}$ and derive the decomposition of the constitutive model, that is, $\varPsi(\bm{\varepsilon}) = \varPsi(\bm{\varepsilon}^\text{t}) + \varPsi(\bm{\varepsilon}^\text{c})$.}
\subsection{Choice of tension-compression asymmetry}\label{Sec:problem-unilateral}
Several existing models for tension-compression asymmetry can be derived in accordance with the variational principle in \eqref{Eq:variational_principle}. Here, we list a few choices:
\begin{itemize}
    \item The original model proposed by Bourdin et al. \cite{Bourdin2000} adopts a symmetric response of the cracked solid, i.e., it assumes that both tension and compression loads contribute equally to cracking. If $\tilde{\mathcal{S}}^\text{t}$ is taken as the set of all symmetric second-order tensors, this model is retrieved.
    \item The volumetric-deviatoric by Amor et al. \cite{Amor2009} assumes that both volumetric expansion and deviatoric deformation contribute to crack propagation but not volumetric compression. $\tilde{\mathcal{S}}^\text{t}$ represents all symmetric second-order tensors with a non-negative trace ($\tr{\tilde{\bm{\varepsilon}}}\geq 0$), i.e. volumetric expansion.
    \item The no-tension by Freddi and Royer-Carfagni \cite{Freddi2010} is a model for masonry-like materials that accounts for the Poisson effect. This model is retrieved when $\tilde{\mathcal{S}}^\text{t}$ is chosen to be the set of all positive semi-definite symmetric tensors. The stress attained at a fully damaged material point is negative semi-definite, which corresponds to the fact that the material does not support tension.
\end{itemize}
    \begin{remark}\label{Re:Intro_Miehe}
        Not all existing approaches that account for the tension-compression asymmetry satisfy the variational principle. For instance, the well known spectral-decomposition model by Miehe et al. \cite{Miehe2010} adopts the elastic energy density split as below, but does not fit into the variational formalism:
        \begin{equation*}
            \begin{aligned}
                \varPsi(\bm{\varepsilon}) &= \varPsi^\text{t}(\bm{\varepsilon}) + \varPsi^\text{c}(\bm{\varepsilon}), \\
                \varPsi^{t}(\bm{\varepsilon}) &= \frac12 \lambda \langle \tr{\bm{\varepsilon}} \rangle_{+} + \mu \bm{\varepsilon}^{t} \cdot \bm{\varepsilon}^{t}, \quad
                \varPsi^{c}(\bm{\varepsilon}) = \frac12 \lambda \langle \tr{\bm{\varepsilon}} \rangle_{-} + \mu \bm{\varepsilon}^{c} \cdot \bm{\varepsilon}^{c},
            \end{aligned}
        \end{equation*}
        where the Macaulay bracket is defined as $\langle\cdot\rangle_{\pm} := (\cdot\pm \lvert\cdot\rvert)/2$, and $\bm{\varepsilon}^{t}$, $\bm{\varepsilon}^{c}$ are obtained by projecting $\bm{\varepsilon}$ onto the space of all symmetric positive/negative semi-definite tensors with respect to the Frobenius norm. Here, contrary to our recent assumption, the constitutive relations cannot be partitioned based on merely $\bm{\varepsilon}^\text{t}$ or $\bm{\varepsilon}^\text{c}$. 
        Consequently, one can not adopt the spectral-decomposition model in the variational formalism. \newline
    \end{remark}
In the following, we provide implementation details of our model based on two choices: (a) {volumetric-deviatoric} and (b) {no-tension}. This means $\tilde{\mathcal{S}}$ can be split into two portions accordingly.
\subsubsection{Volumetric-deviatoric}
We split the transformed strain into volumetric and deviatoric parts,
\begin{subequations}\label{Eq:Amor_split}
    \begin{align}
        {\tilde{\mathcal{S}}}^\text{t}=\Big\{\tilde{\bm{\varepsilon}}\in\tilde{\mathcal{S}}\Big|\tr({\tilde{\bm{\varepsilon}}})\geq 0\Big\}, \\
        {\tilde{\mathcal{S}}}^\text{c}=\Big\{\tilde{\bm{\varepsilon}}\in\tilde{\mathcal{S}}\Big|\tilde{\bm{\varepsilon}}=a\mathbf{1},~a \leq 0\Big\},
    \end{align}
\end{subequations}
where $\mathbf{1}$ is the second-order identity tensor. 
Thus, the strain in the transformed space is decomposed as follows:
\begin{subequations}\label{Eq:Amor_decomp}
    \begin{align}
        \tilde{\boldsymbol{\varepsilon}}^{t} = \frac13\langle\tr{\tilde{\boldsymbol{\varepsilon}}}\rangle_{+}\mathbf{1} 
        + \dev\tilde{\boldsymbol{\varepsilon}},\\
        \tilde{\boldsymbol{\varepsilon}}^{c} = \frac13\langle\tr{\tilde{\boldsymbol{\varepsilon}}}\rangle_{-}\mathbf{1},
    \end{align}
\end{subequations}
where $\dev\tilde{\boldsymbol{\varepsilon}} = \tilde{\boldsymbol{\varepsilon}} - \frac13 (\tr{\tilde{\boldsymbol{\varepsilon}}})\textbf{1}.$
Note that this partitioning satisfies \eqref{Eq:eps_tild_cond}.
\subsubsection*{Decomposition of constitutive model}
The derivatives of $\tilde{\boldsymbol{\varepsilon}}^{t}$ and $\tilde{\boldsymbol{\varepsilon}}^{c}$ with respect to 
$\tilde{\boldsymbol{\varepsilon}}$ define two projection tensors as:
\begin{subequations}\label{Eq:D_tensors}
    \begin{align}
        \mathbb{D}^\text{t} := \frac{\partial\tilde{\bm{\varepsilon}}^\text{t}}{\partial\tilde{\bm{\varepsilon}}} = \frac13 H\bigl[\tr{(\mathbb{C}^{1/2}\boldsymbol{\varepsilon})}\bigr]\mathbf{1}\otimes\mathbf{1}
        + \Bigl(\mathds{1}-\frac13\mathbf{1}\otimes\mathbf{1}\Bigr),\\
        \mathbb{D}^\text{c} := \frac{\partial\tilde{\bm{\varepsilon}}^\text{c}}{\partial\tilde{\bm{\varepsilon}}} = \frac13 H\bigl[-\tr{(\mathbb{C}^{1/2}\boldsymbol{\varepsilon})}\bigr]\mathbf{1}\otimes\mathbf{1},
    \end{align}
\end{subequations}
where $\mathds{1}$ is the fourth-order identity tensor, and $H$ is the Heaviside function such that $H(a) = 1$ if $a\geq 0$, and $H(a)=0$ otherwise.
Then, by applying the chain rule, we obtain:
\begin{equation*}
    \begin{aligned}
        \frac{\partial\bm{\varepsilon}^\text{t}}{\partial\bm{\varepsilon}}=\mathbb{C}^{-1/2}\mathbb{D}^\text{t}\mathbb{C}^{1/2},
        \quad \frac{\partial\bm{\varepsilon}^\text{c}}{\partial\bm{\varepsilon}}=\mathbb{C}^{-1/2}\mathbb{D}^\text{c}\mathbb{C}^{1/2}.
    \end{aligned}
\end{equation*}
Therefore, the decomposition of ${\boldsymbol{\varepsilon}}$ into ${\boldsymbol{\varepsilon}}^{t}~\text{and}~{\boldsymbol{\varepsilon}}^{c}$ reads:
\begin{subequations}\label{Eq:strain_decomp_voldev}
    \begin{align}
        \boldsymbol{\varepsilon}^{t} &= \Bigl(\mathbb{C}^{-1/2}\mathbb{D}^\text{t}\mathbb{C}^{1/2}\Bigr)\boldsymbol{\varepsilon},\\
        \boldsymbol{\varepsilon}^{c} &= \Bigl(\mathbb{C}^{-1/2}\mathbb{D}^\text{c}\mathbb{C}^{1/2}\Bigr)\boldsymbol{\varepsilon}.
    \end{align}
\end{subequations}
Moreover, we derive the stress decomposition as:
\begin{equation}\label{Eq:sigma}
    \begin{aligned}
        \bm{\sigma}^\text{t} := \frac{\partial\varPsi^\text{t}}{\partial\bm{\varepsilon}} = \mathbb{C}\bm{\varepsilon}^\text{t}
        = \mathbb{C}^\text{t}\bm{\varepsilon}, \quad \bm{\sigma}^\text{c} := \frac{\partial\varPsi^\text{c}}{\partial\bm{\varepsilon}} = \mathbb{C}\bm{\varepsilon}^\text{c}
        = \mathbb{C}^\text{c}\bm{\varepsilon}.
    \end{aligned}
\end{equation}
On this basis, we can also split $\mathbb{C}=\frac{\partial^2\varPsi}{\partial\bm{\varepsilon}^2}$ into two portions \[\mathbb{C} = \mathbb{C}^\text{t} + \mathbb{C}^\text{c},\]
where $\mathbb{C}^{t}$ and $\mathbb{C}^\text{c}$ are respectively,
\begin{subequations}\label{Eq:C_decomp}
    \begin{align}
        \mathbb{C}^{t} := \frac{\partial^2 \varPsi^\text{t}}{\partial\bm{\varepsilon}^2} = \bigl(\frac{\partial\bm{\varepsilon}^\text{t}}{\partial{\boldsymbol{\varepsilon}}}\bigr)^\text{T}\mathbb{C}\frac{\partial\bm{\varepsilon}^\text{t}}{\partial{\boldsymbol{\varepsilon}}},\\
        \mathbb{C}^{c} := \frac{\partial^2 \varPsi^\text{c}}{\partial\bm{\varepsilon}^2} = \bigl(\frac{\partial\bm{\varepsilon}^\text{c}}{\partial{\boldsymbol{\varepsilon}}}\bigr)^\text{T}\mathbb{C}\frac{\partial\bm{\varepsilon}^\text{c}}{\partial{\boldsymbol{\varepsilon}}}.
    \end{align}
\end{subequations}
\subsubsection{No-tension}
Taking into account that $\tilde{\bm{\varepsilon}}^\text{t}\cdot\tilde{\bm{\varepsilon}}^\text{c} = 0$, we obtain the following relation for stress-strain in their transformed spaces:
\begin{equation}\label{Eq:mas-isotropy}
    \begin{aligned}
        \tilde{\bm{\sigma}} &:= \frac{\partial\varPsi(\tilde{\bm{\varepsilon}})}{\partial\tilde{\bm{\varepsilon}}} = \frac{\partial}{\partial\tilde{\bm{\varepsilon}}}\Bigl(\frac12\tilde{\bm{\varepsilon}}^\text{t}\cdot\tilde{\bm{\varepsilon}}^\text{t} + \frac12\tilde{\bm{\varepsilon}}^\text{c}\cdot\tilde{\bm{\varepsilon}}^\text{c}\Bigr) = \mathds{1}\tilde{\bm{\varepsilon}}.
    \end{aligned}
\end{equation}
Thus, the stiffness tensor in the transformed spaces is the fourth-order identity tensor, that is, $\tilde{\mathbb{C}}=\frac{\partial\tilde{\bm{\sigma}}}{\partial\tilde{\bm{\varepsilon}}}=\mathds{1}$. On this basis, we can simplify and follow {no-tension} as for the isotropic materials \cite{Freddi2010}. Henceforth, let $(\tilde{\varepsilon}_1, \tilde{\varepsilon}_2, \tilde{\varepsilon}_3)$ be a spectral decomposition of $\tilde{\bm{\varepsilon}}$ such that $\tilde{\varepsilon}_1\geq\tilde{\varepsilon}_2\geq\tilde{\varepsilon}_3$,
\begin{equation}\label{Eq:mas-eps}
    \begin{aligned}
        \tilde{\bm{\varepsilon}} = \overset{3}{\underset{i=1}{\sum}}\tilde{\varepsilon}_i \bm{n}_{(i)} \otimes \bm{n}_{(i)} := \tilde{\varepsilon}_i \bm{M}_i,
    \end{aligned}
\end{equation}
{where $\bm{n}_{(i)}$ represents an eigenvector and parentheses around an index indicate that the usual summation convention is suspended.} 
The strain is decomposed so that $\tilde{\bm{\varepsilon}}^\text{t}$ is a positive definite tensor and is coaxial with $\tilde{\bm{\varepsilon}}$,
\begin{equation}\label{Eq:mas-epspn}
    \begin{aligned}
        \tilde{\bm{\varepsilon}}^{t} = \overset{3}{\underset{i=1}{\sum}} a_i \bm{n}_{(i)} \otimes \bm{n}_{(i)} := a_i \bm{M}_i, \quad \tilde{\bm{\varepsilon}}^{c} = \overset{3}{\underset{i=1}{\sum}} b_i \bm{n}_{(i)} \otimes \bm{n}_{(i)} := b_i \bm{M}_i,
    \end{aligned}
\end{equation}
where $a_i = \langle\tilde{\varepsilon}_i\rangle_{+}$, and $b_i = \tilde{\varepsilon}_i - a_i$. 
\subsubsection*{Decomposition of constitutive model}
A necessary part is to obtain the derivatives of $\tilde{\bm{\varepsilon}}^\text{t}$ and $\tilde{\bm{\varepsilon}}^\text{c}$ with respect to
$\tilde{\bm{\varepsilon}}$.
\begin{equation}\label{Eq:proj-dir-dev}
    \begin{aligned}
        \mathbb{D}^{t} := \frac{\partial\tilde{\bm{\varepsilon}}^t}{\partial{\tilde{\bm{\varepsilon}}}} &= \overset{3}{\underset{i=1}{\sum}} \frac{\partial a_i}{\partial\tilde{\varepsilon}_i} \bm{n}_{(i)} \otimes \bm{n}_{(i)} \otimes\bm{n}_{(i)} \otimes \bm{n}_{(i)} \\
        &+ \frac12\underset{i\neq j}{\sum} \frac{a_j -a_i}{\tilde{\varepsilon}_j - \tilde{\varepsilon}_i} \Biggl(\bm{n}_{(i)} \otimes \bm{n}_{(j)} \otimes \bm{n}_{(i)} \otimes \bm{n}_{(j)} + \bm{n}_{(i)} \otimes \bm{n}_{(j)} \otimes \bm{n}_{(j)} \otimes \bm{n}_{(i)} \Biggr), \\
        \mathbb{D}^{c} := \frac{\partial\tilde{\bm{\varepsilon}}^c}{\partial{\tilde{\bm{\varepsilon}}}} &= \overset{3}{\underset{i=1}{\sum}} \frac{\partial b_i}{\partial\tilde{\varepsilon}_i} \bm{n}_{(i)} \otimes \bm{n}_{(i)} \otimes\bm{n}_{(i)} \otimes \bm{n}_{(i)} \\
        &+ \frac12 \underset{i\neq j}{\sum} \frac{b_j - b_i}{\tilde{\varepsilon}_j - \tilde{\varepsilon}_i} \Biggl(\bm{n}_{(i)} \otimes \bm{n}_{(j)} \otimes \bm{n}_{(i)} \otimes \bm{n}_{(j)} + \bm{n}_{(i)} \otimes \bm{n}_{(j)} \otimes \bm{n}_{(j)} \otimes \bm{n}_{(i)} \Biggr).
    \end{aligned}
\end{equation}
{For sake of brevity, readers are referred to \ref{App} for a closed-form expression of projection tensors in two-dimensional and three-dimensional settings.} \newline
Similar to {volumetric-deviatoric}, the decomposition of $\bm{\varepsilon}$ can be realized as 
\begin{equation}\label{Eq:strain_decomp_masonry}
    \begin{aligned}
        \boldsymbol{\varepsilon}^{t} = \Bigl(\mathbb{C}^{-1/2}\mathbb{D}^\text{t}\mathbb{C}^{1/2}\Bigr)\boldsymbol{\varepsilon}, \quad
        \boldsymbol{\varepsilon}^{c} = \Bigl(\mathbb{C}^{-1/2}\mathbb{D}^\text{c}\mathbb{C}^{1/2}\Bigr)\boldsymbol{\varepsilon},
    \end{aligned}
\end{equation}
with $\mathbb{D}^\text{t}, \mathbb{D}^\text{c}$ as in \eqref{Eq:proj-dir-dev}. \newline
In addition, we write the stress decomposition as
\begin{equation}\label{Eq:sigma_masonry}
    \begin{aligned}
        \sigma^\text{t}_{mn} := \frac{\partial\varPsi^\text{t}}{\partial\varepsilon_{mn}} = \frac{\partial\varPsi^\text{t}}{\partial\tilde{\varepsilon}_{ij}^\text{t}}\frac{\partial\tilde{\varepsilon}_{ij}^\text{t}}{\partial\tilde{\varepsilon}_{kl}}\frac{\partial\tilde{\varepsilon}_{kl}}{\partial\varepsilon_{mn}}=\tilde{\varepsilon}_{ij}^\text{t} D^\text{t}_{ijkl}\mathbb{C}^{1/2}_{klmn} = \tilde{\varepsilon}^\text{t}_{kl} C^{1/2}_{klmn} = C^{1/2}_{mnkl}C^{1/2}_{klpq}{\varepsilon}^\text{t}_{pq}, \\
        \sigma^\text{c}_{mn} := \frac{\partial\varPsi^\text{c}}{\partial\varepsilon_{mn}} = \frac{\partial\varPsi^\text{c}}{\partial\tilde{\varepsilon}_{ij}^\text{c}}\frac{\partial\tilde{\varepsilon}_{ij}^\text{c}}{\partial\tilde{\varepsilon}_{kl}}\frac{\partial\tilde{\varepsilon}_{kl}}{\partial\varepsilon_{mn}}=\tilde{\varepsilon}_{ij}^\text{c} D^\text{c}_{ijkl}\mathbb{C}^{1/2}_{klmn} = \tilde{\varepsilon}^\text{c}_{kl} C^{1/2}_{klmn} = C^{1/2}_{mnkl}C^{1/2}_{klpq}{\varepsilon}^\text{c}_{pq},
    \end{aligned}
\end{equation}
or in tensorial form:
\begin{equation}\label{Eq:sigma_masonry_tensor}
    \begin{aligned}
        \bm{\sigma}^\text{t} := \frac{\partial\varPsi^\text{t}}{\partial\bm{\varepsilon}} = \mathbb{C} \bm{\varepsilon}^\text{t}, \quad \bm{\sigma}^\text{c} := \frac{\partial\varPsi^\text{c}}{\partial\bm{\varepsilon}} = \mathbb{C}\bm{\varepsilon}^\text{c}.
    \end{aligned}
\end{equation}
To split $\mathbb{C}=\frac{\partial^2\varPsi}{\partial\bm{\varepsilon}^2}$ into two portions, we take the derivatives of \eqref{Eq:sigma_masonry} as:
\begin{equation}\label{Eq:C_masonry}
    \begin{aligned}
        \frac{\partial^2\varPsi^\text{t}}{\partial\varepsilon_{pq}\partial\varepsilon_{mn}} = \frac{\partial}{\partial\varepsilon_{pq}}\bigl(\tilde{\varepsilon}^\text{t}_{ij}D^\text{t}_{ijkl}\bigr)C^{1/2}_{klmn} = C^{1/2}_{mnkl}\frac{\partial}{\partial\tilde{\varepsilon}_{rs}}\bigl(\tilde{\varepsilon}^\text{t}_{ij}D^\text{t}_{ijkl}\bigr)\frac{\partial\tilde{\varepsilon}_{rs}}{\partial\varepsilon_{pq}} 
        = C^{1/2}_{mnkl}\frac{\partial\tilde{\varepsilon}^\text{t}_{kl}}{\partial\tilde{\varepsilon}_{rs}}\frac{\partial\tilde{\varepsilon}_{rs}}{\partial\varepsilon_{pq}} = C^{1/2}_{mnkl} D^\text{t}_{klrs}C^{1/2}_{rspq}, \\
        \frac{\partial^2\varPsi^\text{c}}{\partial\varepsilon_{pq}\partial\varepsilon_{mn}} = \frac{\partial}{\partial\varepsilon_{pq}}\bigl(\tilde{\varepsilon}^\text{c}_{ij}D^\text{c}_{ijkl}\bigr)C^{1/2}_{klmn} = C^{1/2}_{mnkl}\frac{\partial}{\partial\tilde{\varepsilon}_{rs}}\bigl(\tilde{\varepsilon}^\text{c}_{ij}D^\text{c}_{ijkl}\bigr)\frac{\partial\tilde{\varepsilon}_{rs}}{\partial\varepsilon_{pq}} 
        = C^{1/2}_{mnkl}\frac{\partial\tilde{\varepsilon}^\text{c}_{kl}}{\partial\tilde{\varepsilon}_{rs}}\frac{\partial\tilde{\varepsilon}_{rs}}{\partial\varepsilon_{pq}} = C^{1/2}_{mnkl} D^\text{c}_{klrs}C^{1/2}_{rspq},
    \end{aligned}
\end{equation}
or in tensorial form:
\begin{equation}\label{Eq:C_masonry_tensor}
    \begin{aligned}
        \mathbb{C}^\text{t} = \mathbb{C}^{1/2} \mathbb{D}^\text{t} \mathbb{C}^{1/2}, \quad \mathbb{C}^\text{c} = \mathbb{C}^{1/2} \mathbb{D}^\text{c} \mathbb{C}^{1/2}.
    \end{aligned}
\end{equation}
\begin{remark}\label{Re:tensor_operation}
Note that in \eqref{Eq:sigma_masonry}, $\tilde{\varepsilon}^\text{t}_{ij}D^\text{t}_{ijkl} = \tilde{\varepsilon}^\text{t}_{kl}$ and $\tilde{\varepsilon}^\text{c}_{ij}D^\text{c}_{ijkl} = \tilde{\varepsilon}^\text{c}_{kl}$. This is achieved because of $\frac{\partial a_i}{\partial\tilde{\varepsilon}_i} = H(\tilde{\varepsilon_i})$ and $\frac{\partial b_i}{\partial\tilde{\varepsilon}_i} = H(-\tilde{\varepsilon_i})$ in simplifying \eqref{Eq:proj-dir-dev}. 
Also note in \eqref{Eq:C_masonry} to avoid running into some higher-order tensors, 
rather than employing the product rule, we first compute the aforementioned contractions
and thereafter we proceed with taking their derivatives, that is,
\begin{equation*}
    \begin{aligned}
        \frac{\partial}{\partial\tilde{\varepsilon}_{rs}}\bigl(\tilde{\varepsilon}^\text{t}_{ij}D^\text{t}_{ijkl}\bigr) = \frac{\partial\tilde{\varepsilon}^\text{t}_{kl}}{\partial\tilde{\varepsilon}_{rs}}, \quad
        \frac{\partial}{\partial\tilde{\varepsilon}_{rs}}\bigl(\tilde{\varepsilon}^\text{c}_{ij}D^\text{t}_{ijkl}\bigr) = \frac{\partial\tilde{\varepsilon}^\text{c}_{kl}}{\partial\tilde{\varepsilon}_{rs}}.
    \end{aligned}
\end{equation*}
\end{remark}
\subsection{Type of anisotropy / material symmetry}\label{Sec:problem-anisotropy}
The presented method is general enough to apply for an arbitrary type of anisotropic elasticity tensors. Here we present three types of anisotropy, (a) cubic symmetry, (b) orthotropy and (c) full anisotropy, in Kelvin-matrix notation ($\undertilde{\bm{C}}$). 
Note that the components of $\undertilde{\bm{C}}$ are obtained assuming the orthonormal basis vectors are orientated along the material principal axes. Otherwise, a rotation of $\undertilde{\bm{C}}$ is required, see Remark \ref{Re:rotation}.
\paragraph{Cubic symmetry}
\begin{equation*}
    \begin{aligned}
        \undertilde{\bm{C}} = 
        \begin{bmatrix}
            \mathbb{C}_{1111} & \mathbb{C}_{1122} & \mathbb{C}_{1133} & 0 & 0 & 0 \\
            & \mathbb{C}_{1111} & \mathbb{C}_{2233} & 0 & 0 & 0 \\
            &  & \mathbb{C}_{1111} & 0 & 0 & 0 \\               &  &  & {2}\mathbb{C}_{2323} & 0 & 0 \\
            &  &  &  & {2}\mathbb{C}_{1313} & 0 \\
            &  &  &  &  & {2}\mathbb{C}_{1212}
        \end{bmatrix}
    \end{aligned}
\end{equation*}
\paragraph{Orthotropy}
    \begin{equation*}
        \begin{aligned}
            \undertilde{\bm{C}} = 
            \begin{bmatrix}
                \mathbb{C}_{1111} & \mathbb{C}_{1122} & \mathbb{C}_{1133} & 0 & 0 & 0 \\
                 & \mathbb{C}_{2222} & \mathbb{C}_{2233} & 0 & 0 & 0 \\
                 &  & \mathbb{C}_{3333} & 0 & 0 & 0 \\
                 &  &  & {2}\mathbb{C}_{2323} & 0 & 0 \\
                 &  &  &  & {2}\mathbb{C}_{1313} & 0 \\
                 &  &  &  &  & {2}\mathbb{C}_{1212}
            \end{bmatrix}
        \end{aligned}
    \end{equation*}
\paragraph{Full anisotropy}
\begin{equation*}
    \begin{aligned}
        \undertilde{\bm{C}} = 
        \begin{bmatrix}
            \mathbb{C}_{1111} & \mathbb{C}_{1122} & \mathbb{C}_{1133} & \sqrt{2}\mathbb{C}_{1123} & \sqrt{2}\mathbb{C}_{1113} & \sqrt{2}\mathbb{C}_{1112} \\
             & \mathbb{C}_{2222} & \mathbb{C}_{2233} & \sqrt{2}\mathbb{C}_{2223} & \sqrt{2}\mathbb{C}_{2213} & \sqrt{2}\mathbb{C}_{2212} \\
             &  & \mathbb{C}_{3333} & \sqrt{2}\mathbb{C}_{3323} & \sqrt{2}\mathbb{C}_{3313} & \sqrt{2}\mathbb{C}_{3312} \\
             &  &  & {2}\mathbb{C}_{2323} & 2\mathbb{C}_{2313} & 2\mathbb{C}_{2312} \\
             &  &  &  & {2}\mathbb{C}_{1313} & 2\mathbb{C}_{1312} \\
             &  &  &  &  & {2}\mathbb{C}_{1212}
        \end{bmatrix}
    \end{aligned}
\end{equation*}
\begin{remark}\label{Re:rotation}
We summarize how to transform $\undertilde{\bm{C}}$ to a coordinate system that does not coincide with the material principal axes. 
Below, $\undertilde{\mathbf{C}}^{\prime}$ is obtained via a transformation to a new basis ($x^\prime, y^\prime, z^\prime$):
\begin{equation}\label{Eq:rotation-p}
    \begin{aligned}
        \undertilde{\mathbf{C}}^{\prime} = \mathbf{P}\undertilde{\mathbf{C}}\mathbf{P}^\text{T},
    \end{aligned}
\end{equation}
where $\mathbf{P}$ denotes the transformation matrix.
In a special case, when the rotation is around the $z$-axis by an angle $\alpha$, $\mathbf{P}$ reduces to:
\begin{equation*}
    \begin{aligned}
        \mathbf{P} =
        \begin{bmatrix}
            c^2 & s^2 & 0 & 0 & 0 & \sqrt{2}cs \\
            s^2 & c^2 & 0 & 0 & 0 & -\sqrt{2}cs \\
            0 & 0 & 1 & 0 & 0 & 0 \\
            0 & 0 & 0 & c & -s & 0 \\
            0 & 0 & 0 & s & c & 0 \\
            -\sqrt{2}cs & \sqrt{2}cs & 0 & 0 & 0 & c^2-s^2 \\
        \end{bmatrix},
    \end{aligned}
\end{equation*}
where $c=\cos{\alpha}$, and $s=\sin{\alpha}$.
\end{remark}
    \section{Numerical examples}\label{Sec:numerical-examples}
\begin{figure}[htbp]
	\centering
	\subfigure[tensile test]{
		\begin{minipage}{0.45\linewidth}
			\centering
			\includegraphics[width=2.5in]{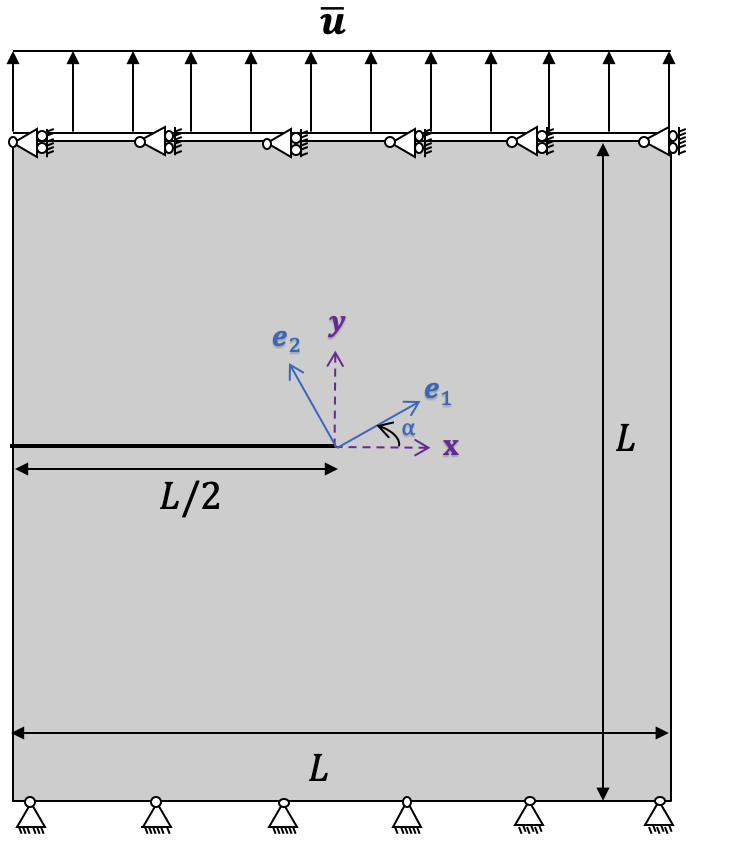}
		\end{minipage}\label{Fig:geom_tensile}
	}
	\subfigure[shear test]{
		\begin{minipage}{0.45\linewidth}
			\centering
			\includegraphics[width=2.5in]{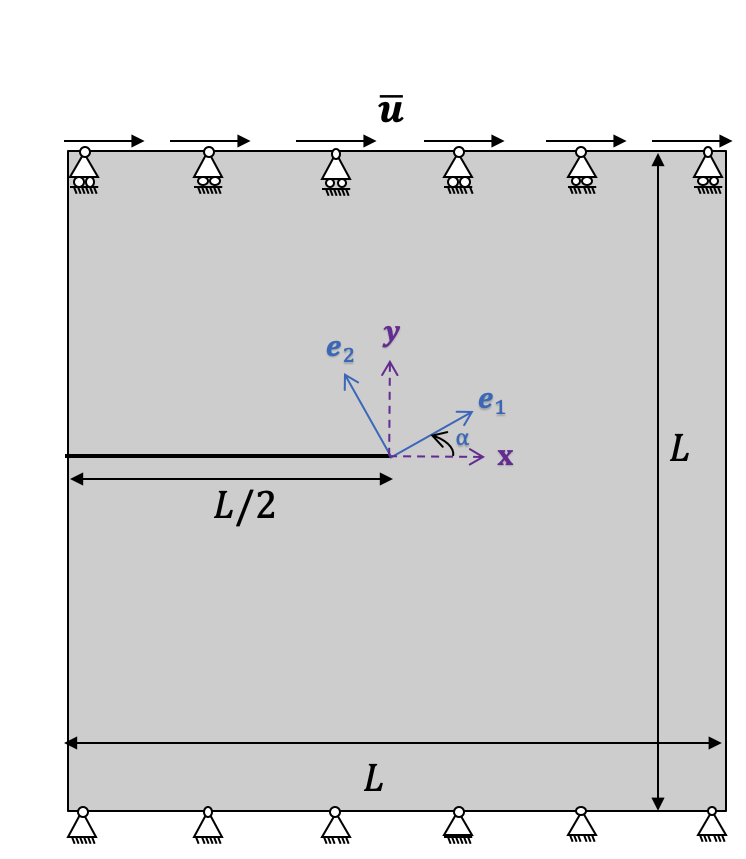}
		\end{minipage}\label{Fig:geom_shear}
	}
	\caption{Schematic of two cracked square plate under a displacement (a) tensile and (b) shear loading. An  incremental displacement loading is applied on the top edges. The material principal axes ($\bm{e}_1,\bm{e}_2$) have an angle of $\alpha$ with the original coordinate system ($\bm{x},\bm{y}$).}
	\label{Fig:Configuration}
\end{figure}

In order to verify our model, here we present two sets of numerical examples, a single edge notched tensile test and a shear test\footnote{The present model is implemented in the phase field process of an open source code, OpenGeoSys~\cite{Yoshioka2019}. Further information on the code and simulation examples are freely accessible at \url{https://www.opengeosys.org/} and~\cite{Bilke2019337}.}. 
The geometry and boundary conditions for both examples are shown in Figure \ref{Fig:geom_tensile} and \ref{Fig:geom_shear}, respectively. \newline
We consider a square plate ($0.001\,\text{m} \times 0.001\, \text{m}$) with an initial horizontal crack placed at the middle height from the left outer surface to the center of the specimen. The specimen was discretized by an unstructured mesh with approximately 41K and 64K standard {triangular} elements, respectively, for the tensile and shear tests. For each discretization, the mesh size of $h=4\times 10^{-6}$ m was used for the critical region of the expected crack path. The length scale parameter was chosen as {$\ell = 8\times 10^{-6}$} m which is small enough with respect to the specimen dimensions. Plane strain condition was assumed. Finally, the \ATtwo{} model was used in all simulations. \newline
The boundary conditions are as follows. The displacement along any direction on the bottom edge ($y=-L/2$) was fixed to zero. 
Also, the displacement at the top edge ($y = L/2$) was prescribed along the $y$-direction for the tensile test and the $x$-direction for the shear test, where the other direction was taken to be zero for both loading setup. 
The computation was performed with a monotonic displacement controlled loading with a constant displacement increment $\Delta\bm{u} = 2\times 10^{-7}$ m until $\Bar{\bm{u}} = 2 \times 10^{-6}$ m. Thereafter, the increment was adjusted to $\Delta\bm{u} = 1 \times 10^{-8}$ m until the last step at $\Bar{\bm{u}} = 2 \times 10^{-5}$ m. Also, the parameters used for an orthotropic material are listed in Table \ref{Tabs:material_properties}. Note that in both examples the critical surface energy $G_c$ was taken to be independent of the material orientation.
\begin{table}[ht]
    \centering
    \begin{tabular}{c c c c}
    \hline
        Material Properties & Value & Unit \\ \hline
        $\nu_{23},\nu_{13},\nu_{12}$ & 0.17, 0.3, 0.52 & - \\ \hline
        $E_1,E_2,E_3$ & 210, 70, 210 & GPa \\ \hline
        $G_{23}, G_{13}, G_{12}$ & 46.63, 80.77, 46.63 & GPa \\ \hline
        $G_c$ & 2700 & N/m \\ \hline
\end{tabular}
    \caption{Orthotropic material properties used in this study, which were taken from van Dijk et al \cite{DIJK2020}.}
    \label{Tabs:material_properties}
\end{table}
\subsection{Tensile test}
\begin{figure}[htbp]
\centering %
\subfigure[\texttt{volumetric-deviatoric}, $\alpha = -\pi/4$]{\includegraphics[width=0.25\linewidth]{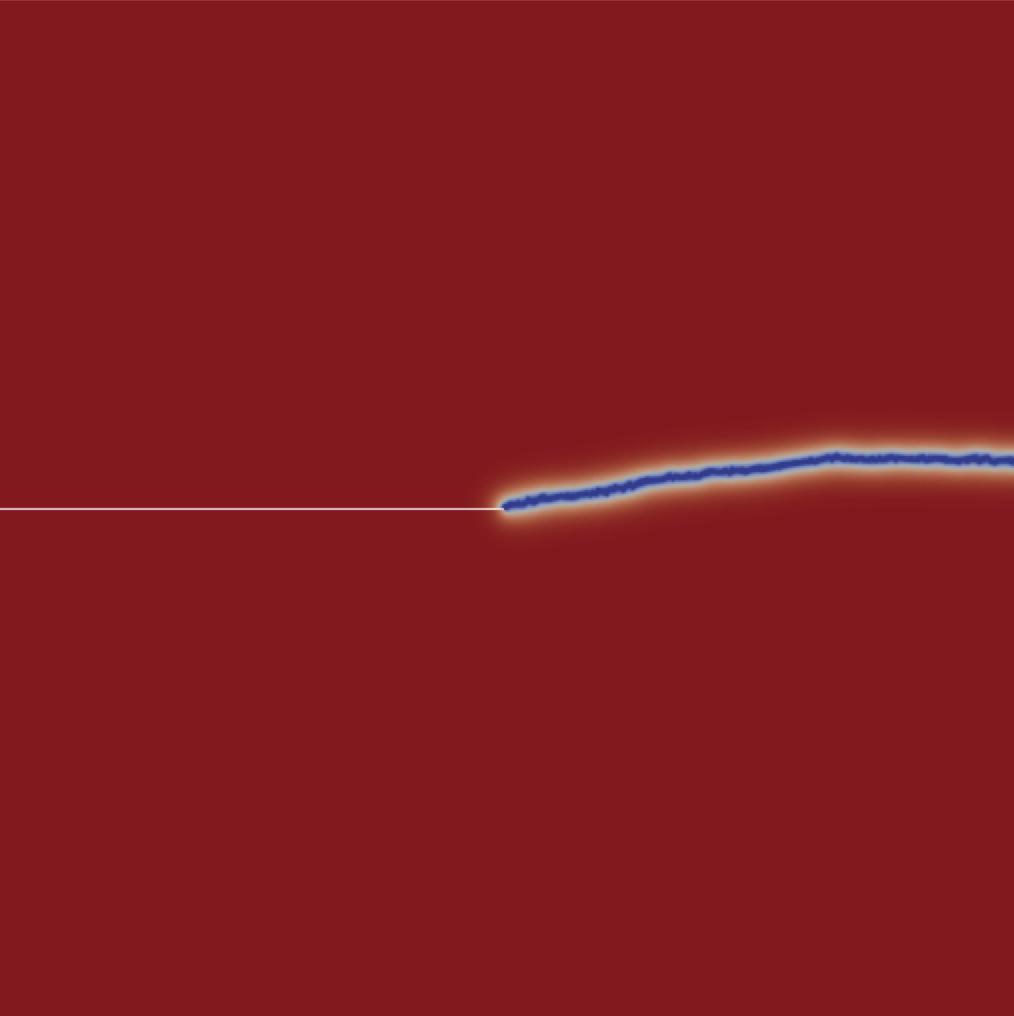}}
\quad
\subfigure[\texttt{no-tension}, $\alpha = -\pi/4$]{\includegraphics[width=0.25\linewidth]{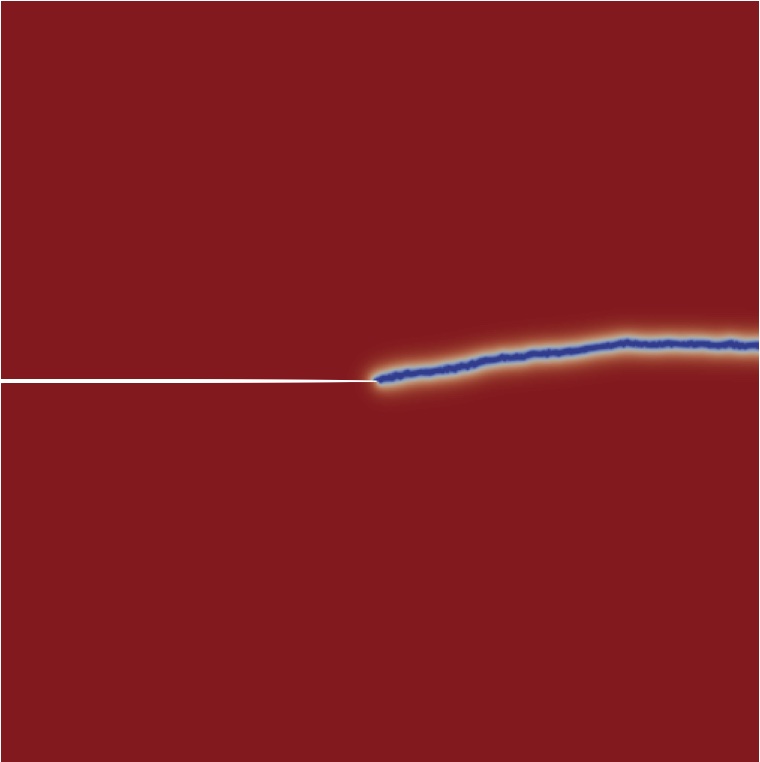}}
\vskip\baselineskip
\subfigure[\texttt{volumetric-deviatoric}, $\alpha = 0$]{\includegraphics[width=0.25\linewidth]{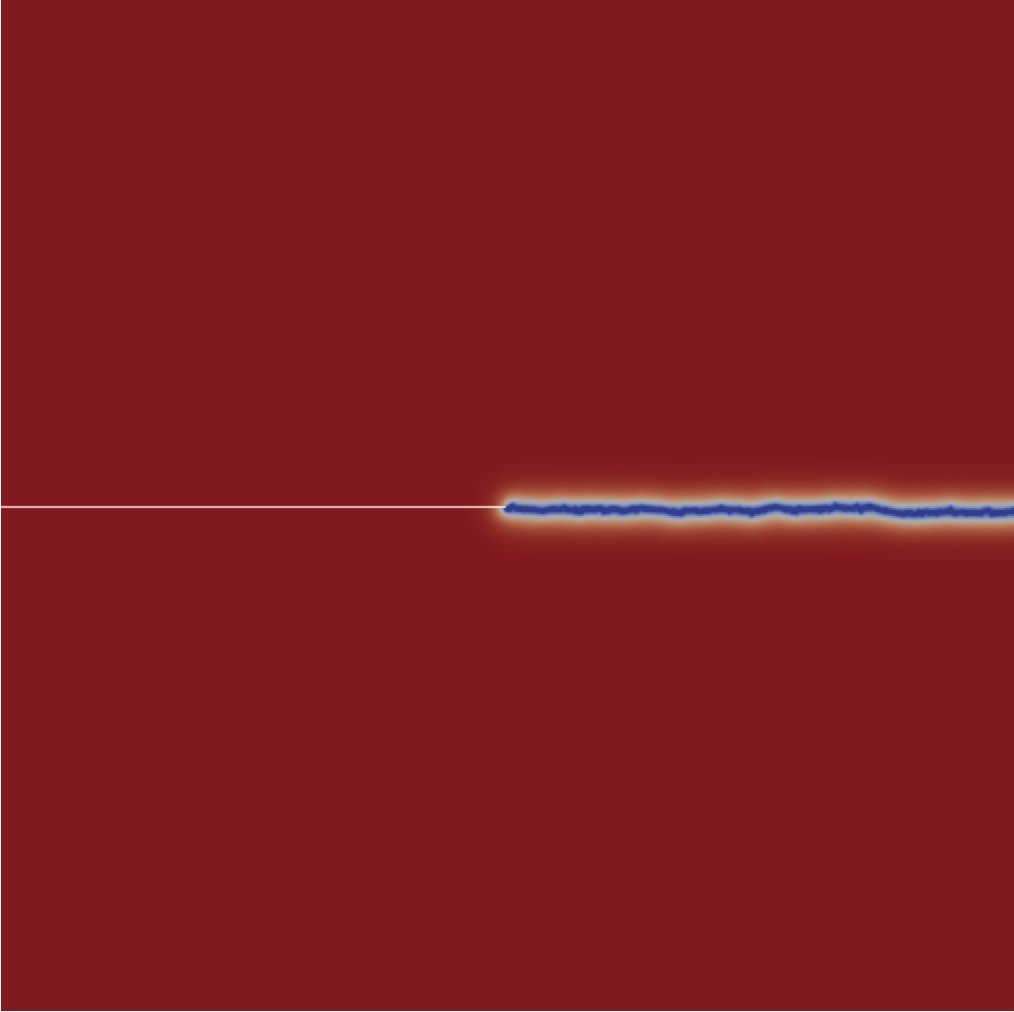}}
\quad
\subfigure[\texttt{no-tension}, $\alpha = 0$]{\includegraphics[width=0.25\linewidth]{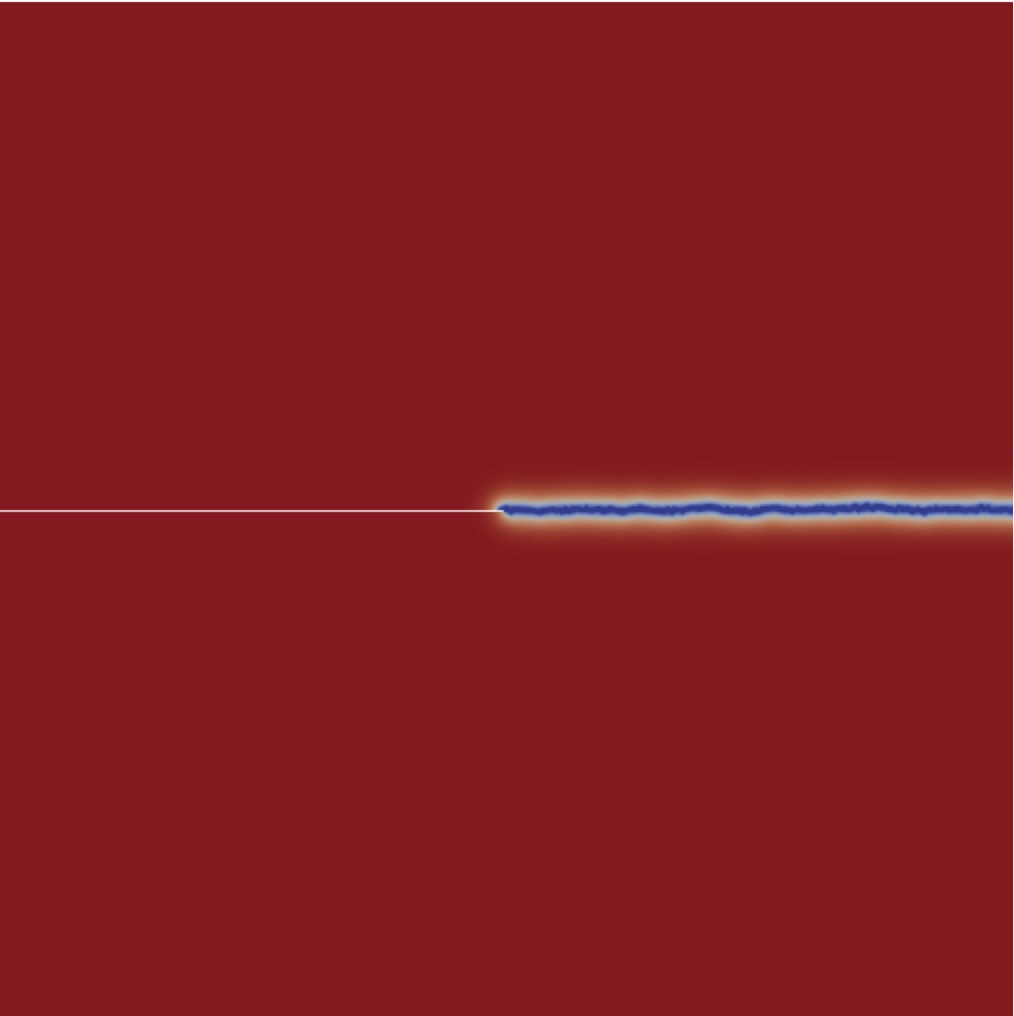}}
\vskip\baselineskip
\subfigure[\texttt{volumetric-deviatoric}, $\alpha = \pi/4$]{\includegraphics[width=0.25\linewidth]{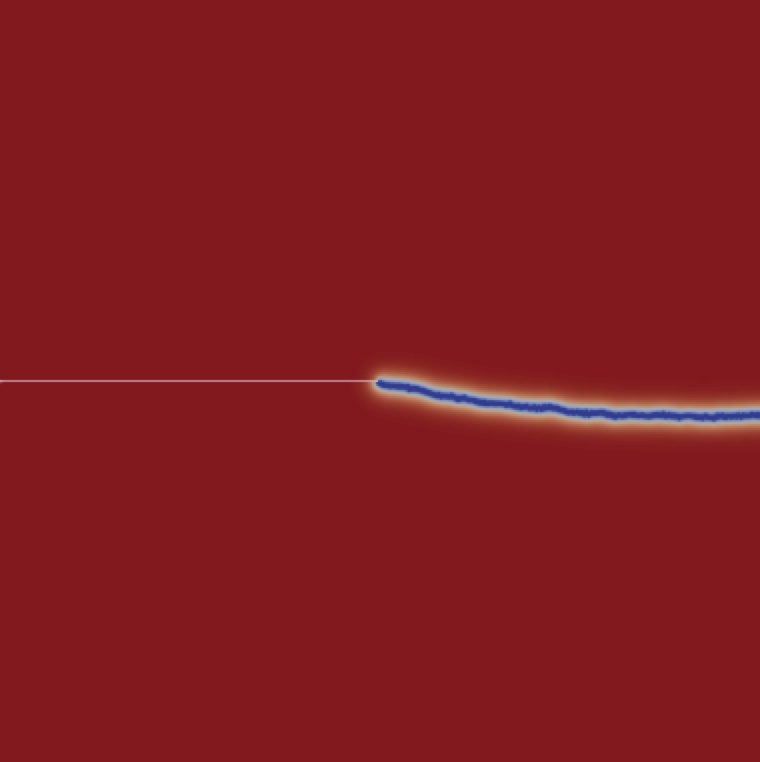}}
\quad
\subfigure[\texttt{no-tension}, $\alpha = \pi/4$]{\includegraphics[width=0.25\linewidth]{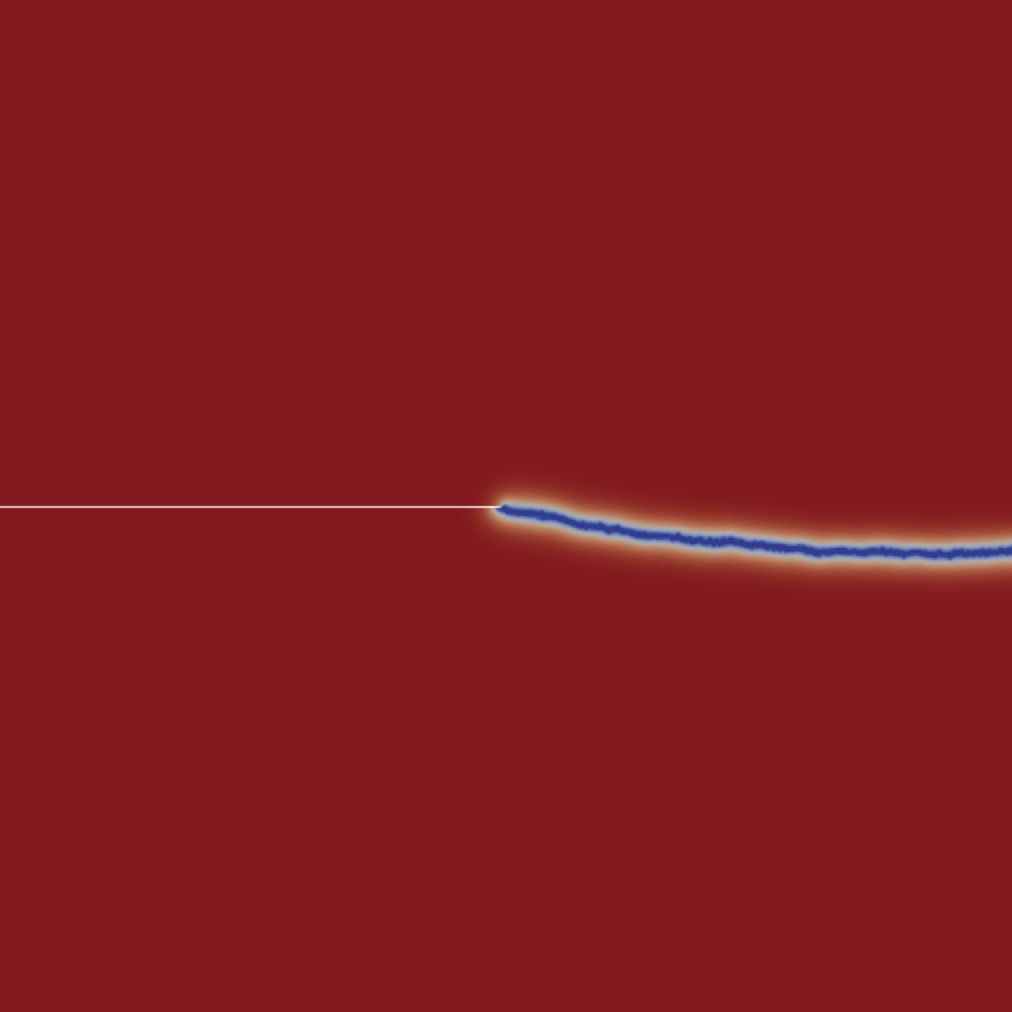}}
\vskip\baselineskip
\subfigure[\texttt{volumetric-deviatoric}, $\alpha = \pi/2$]{\includegraphics[width=0.25\linewidth]{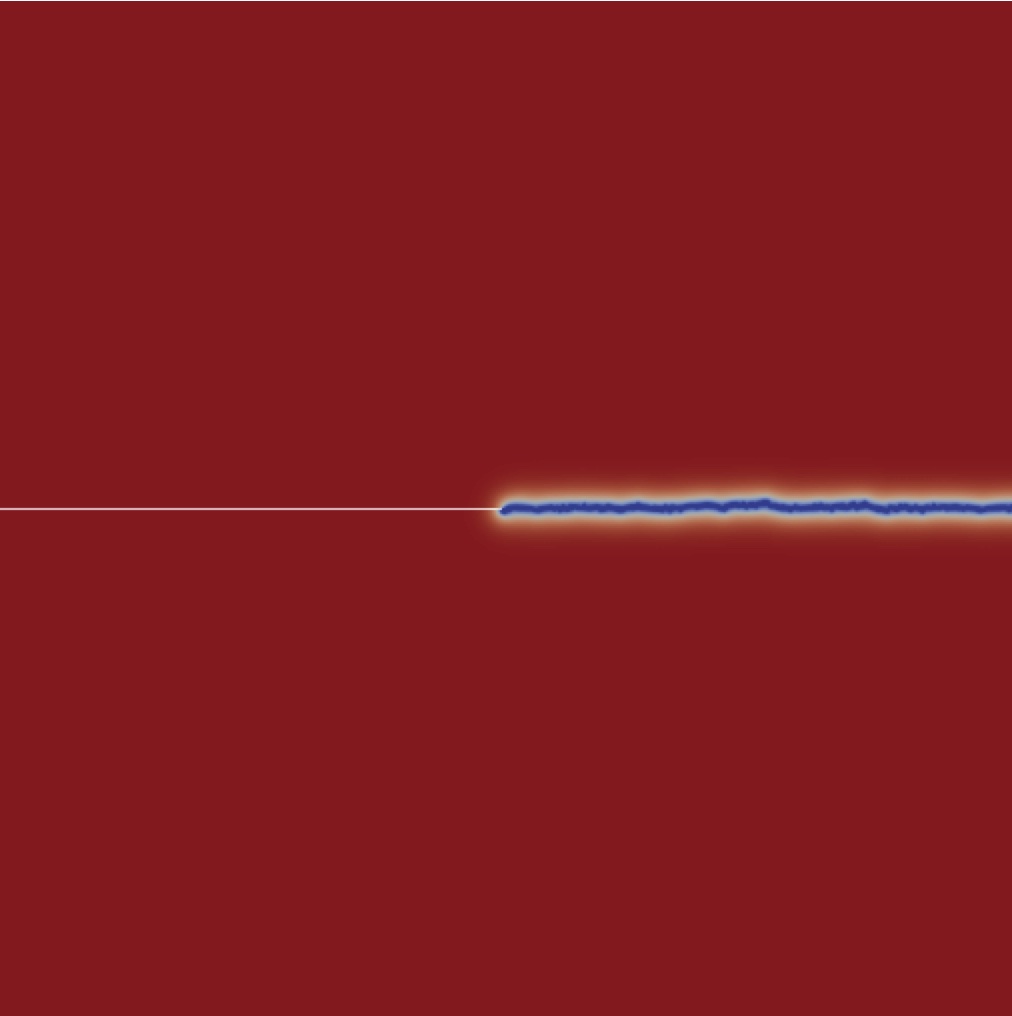}}
\quad
\subfigure[\texttt{no-tension}, $\alpha = \pi/2$]{\includegraphics[width=0.25\linewidth]{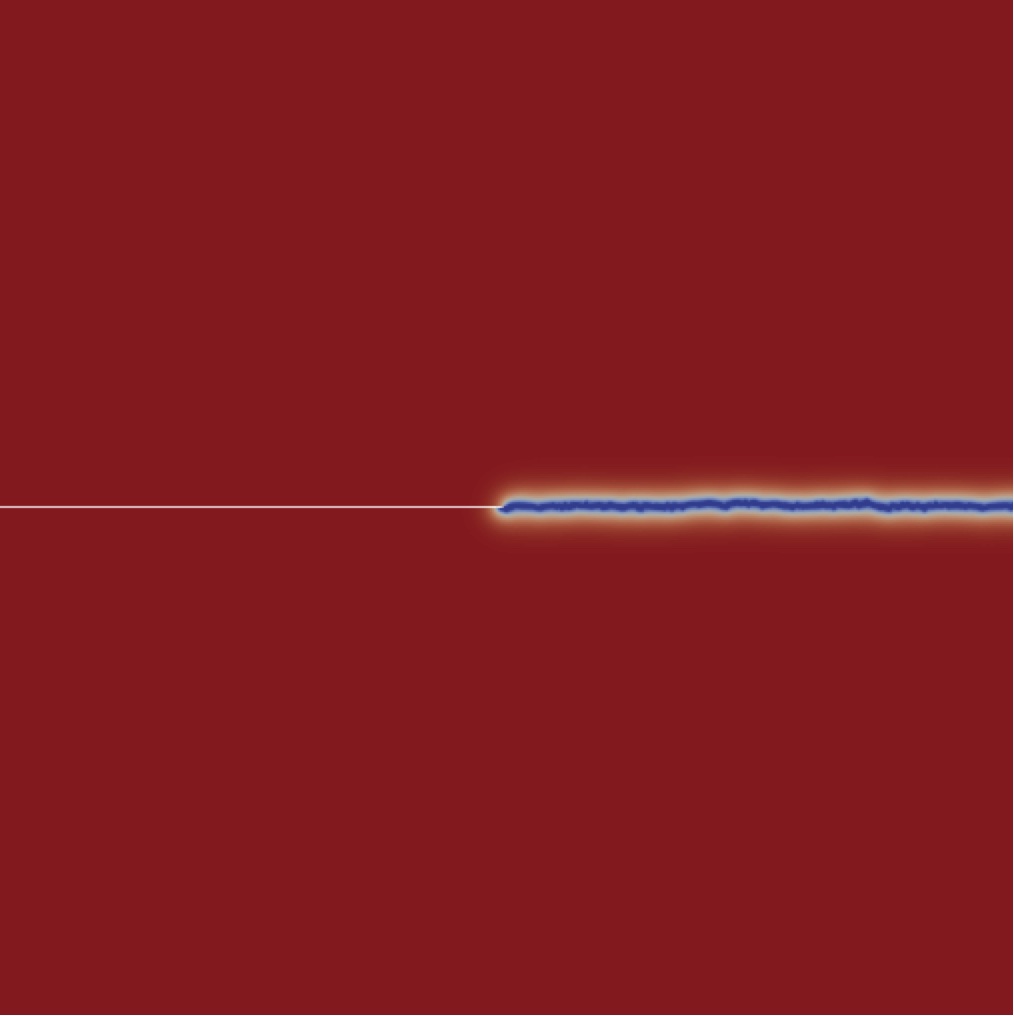}}
\caption{Cracked square plate under tension. Phase field contours of orthotropic {volumetric-deviatoric} and {no-tension} at the last step with with $\Bar{\bm{u}} = 2\times 10^{-5}$ m. The same parameters were used for both models. Under tensile loading, we observe a similar crack response for both models. With four material orientation angles tested for both models, the crack propagates along the horizontal axis for $\alpha=0$ and $\alpha=\pi/2$. For $\alpha=\pm\pi/4$, the crack deviates slightly for both models as a consequence of orthotropy.}
\label{fig:tension-pf}
\end{figure}

\begin{figure}[htbp]
\centering %
\subfigure[$\alpha = -\pi/4$, $\bm{u} \approx 0.37\Bar{\bm{u}}$]{\includegraphics[width=0.4\linewidth]{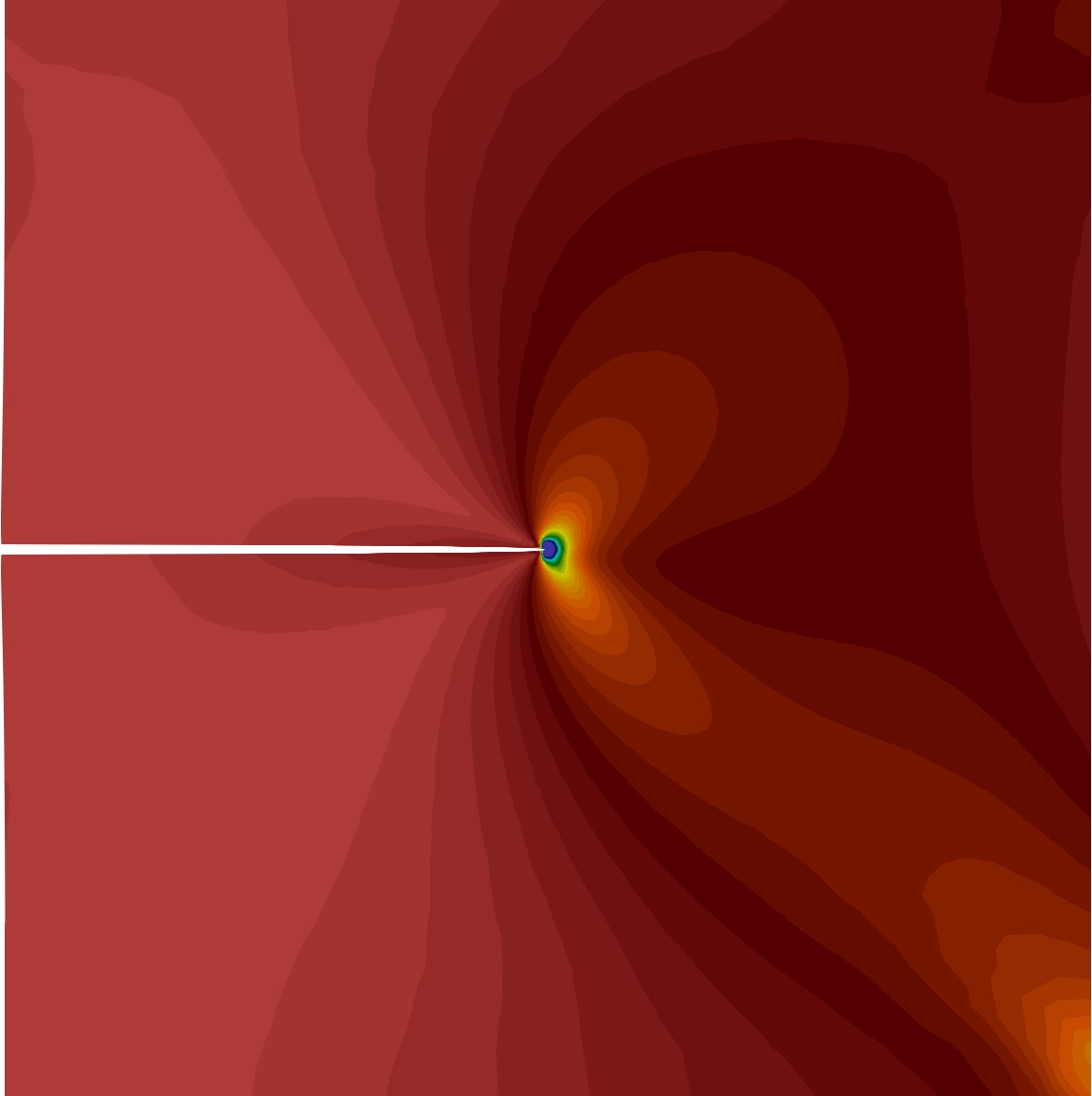}}
\quad
\subfigure[$\alpha = 0$, $\bm{u} \approx 0.47\Bar{\bm{u}}$]{\includegraphics[width=0.4\linewidth]{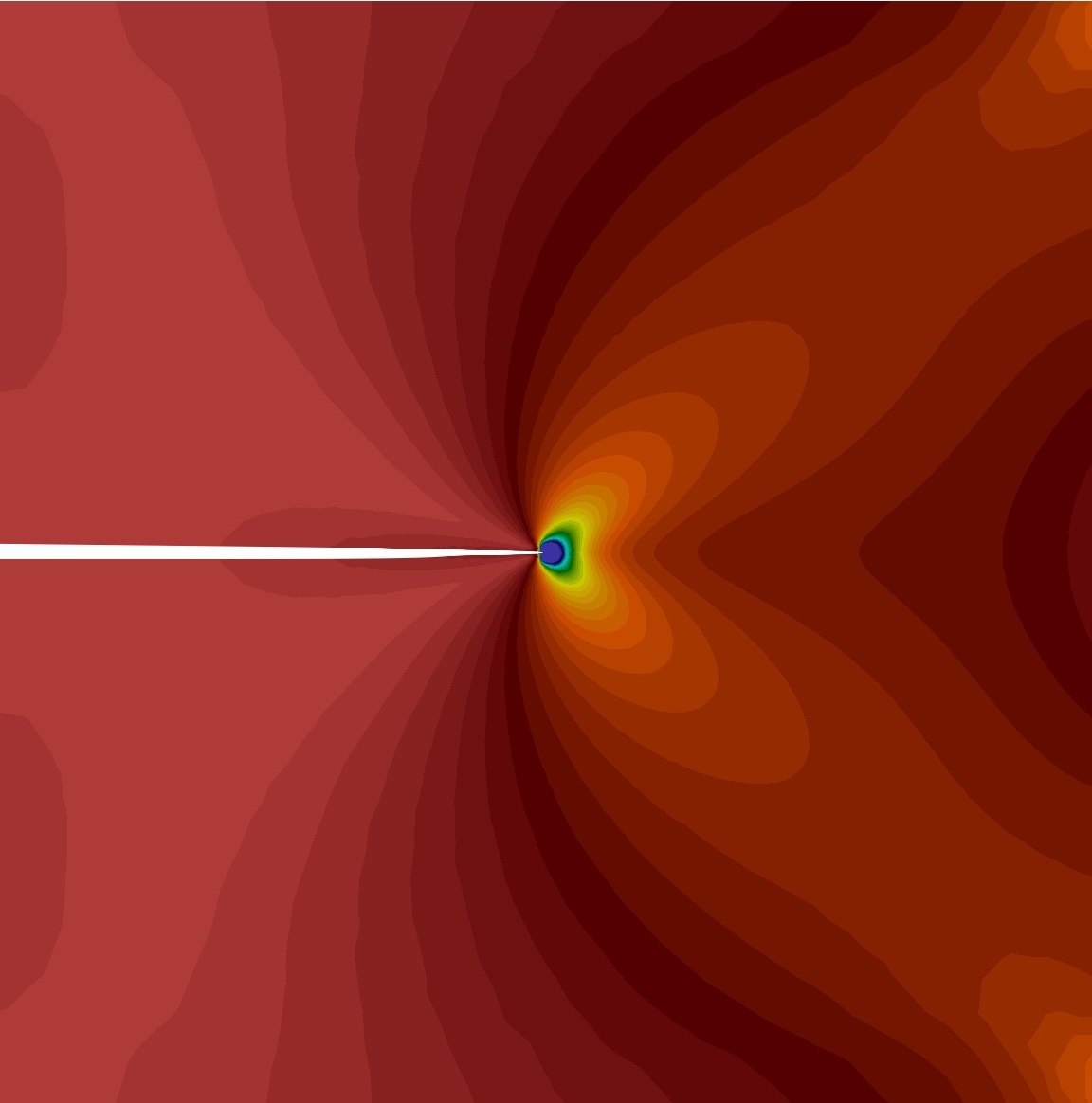}}
\quad
\subfigure[$\alpha = \pi/4$, $\bm{u} \approx 0.36 \Bar{\bm{u}}$]{\includegraphics[width=0.4\linewidth]{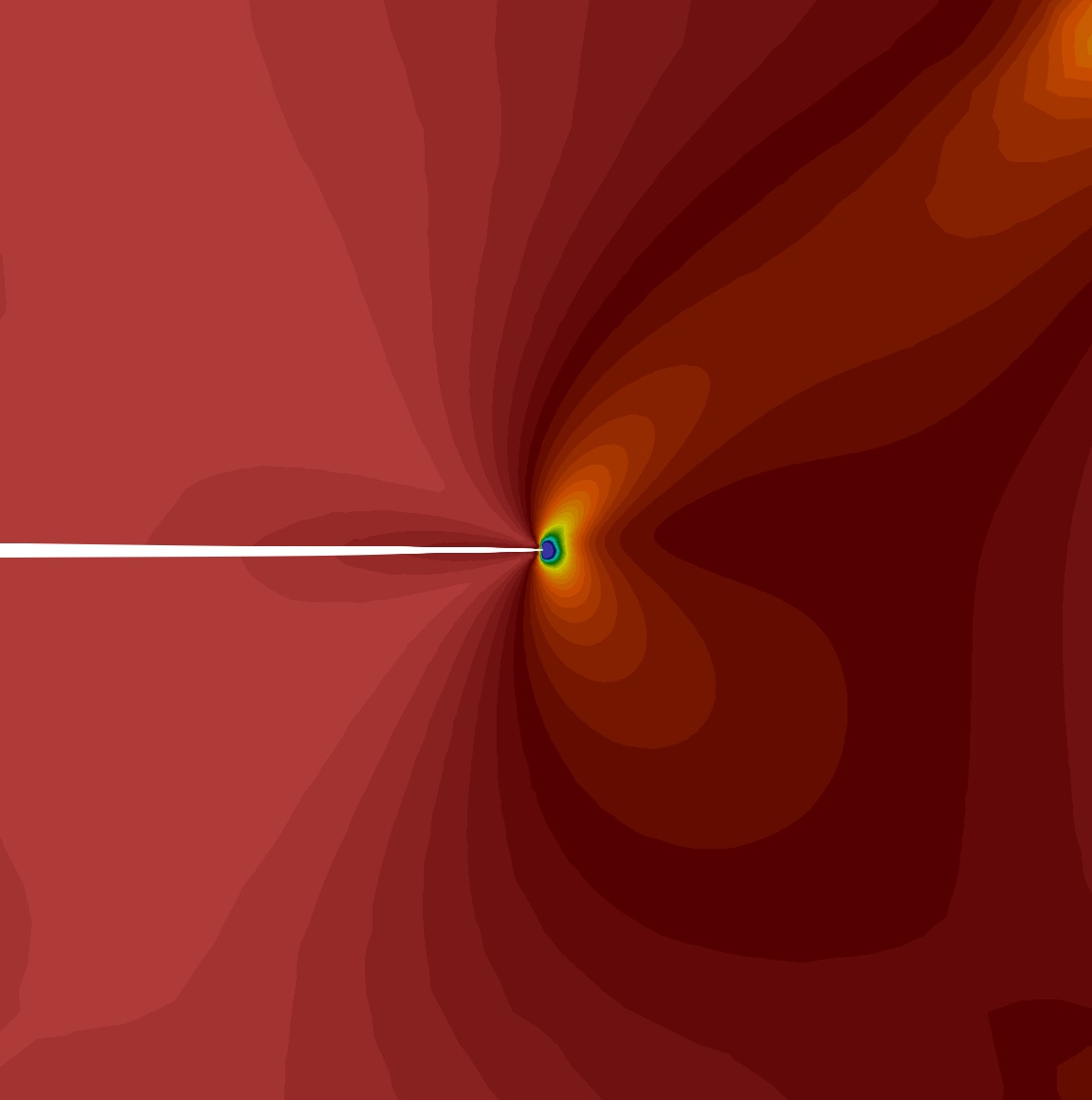}}
\quad
\subfigure[$\alpha = \pi/2$, $\bm{u} \approx 0.30 \Bar{\bm{u}}$]{\includegraphics[width=0.4\linewidth]{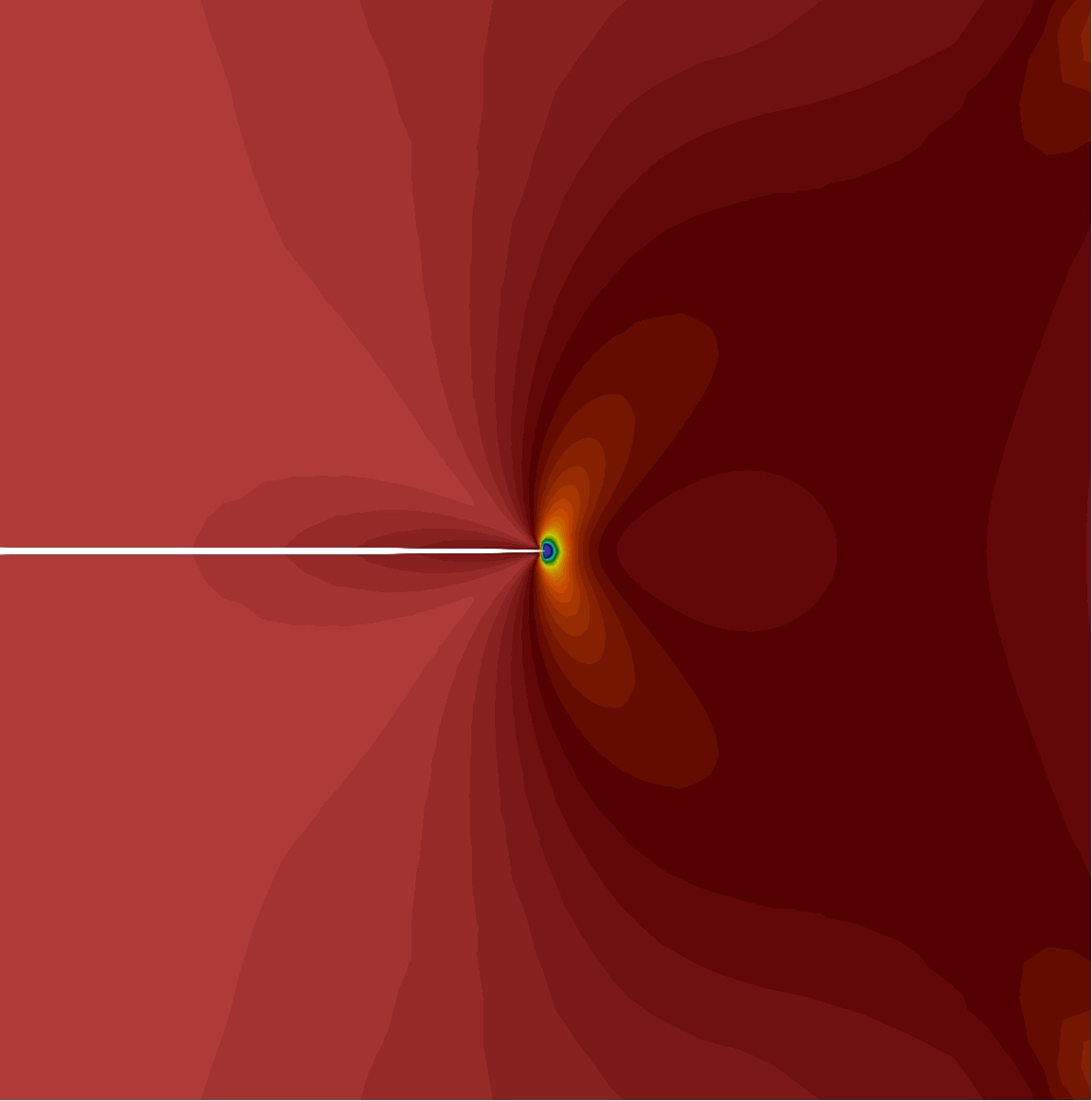}}\\
\subfigure{\includegraphics[width=0.25\linewidth]{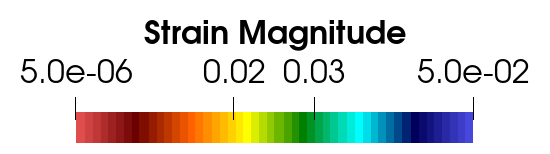}}
\caption{The magnitude of the strain tensor at the initiation of the crack for orthotropic {volumetric-deviatoric} with four material orientation angles. A symmetric behavior is observed for $\alpha=0$ and $\alpha=\pi/2$, while the magnitude is higher for $\alpha=0$ as the material has more stiffness along the horizontal axis. For $\alpha=\pm\pi/4$, asymmetry is observed where the magnitude is higher along the material principal axis with more stiffness.}
\label{fig:tension-eps-ovd}
\end{figure}

\begin{figure}[htbp]
\centering %
\subfigure[\texttt{volumetric-deviatoric}, $\alpha = -\pi/4$]{\includegraphics[width=0.35\linewidth]{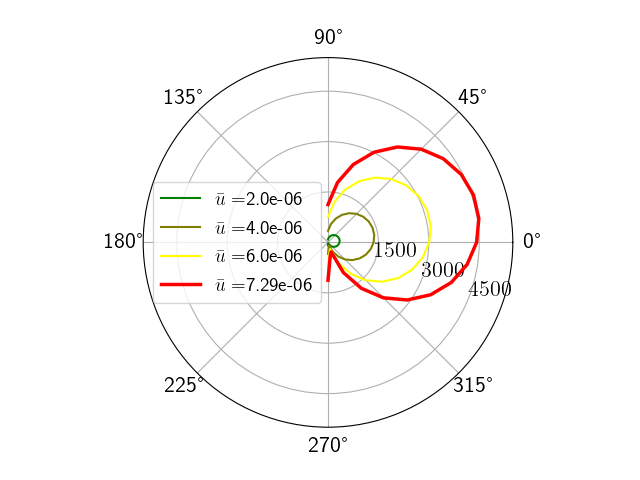}}
\quad
\subfigure[\texttt{no-tension}, $\alpha = -\pi/4$]{\includegraphics[width=0.35\linewidth]{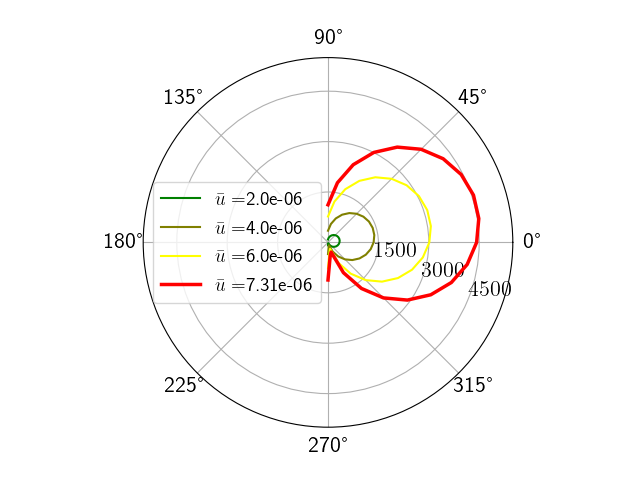}}
\vskip\baselineskip
\subfigure[\texttt{volumetric-deviatoric}, $\alpha = 0$]{\includegraphics[width=0.35\linewidth]{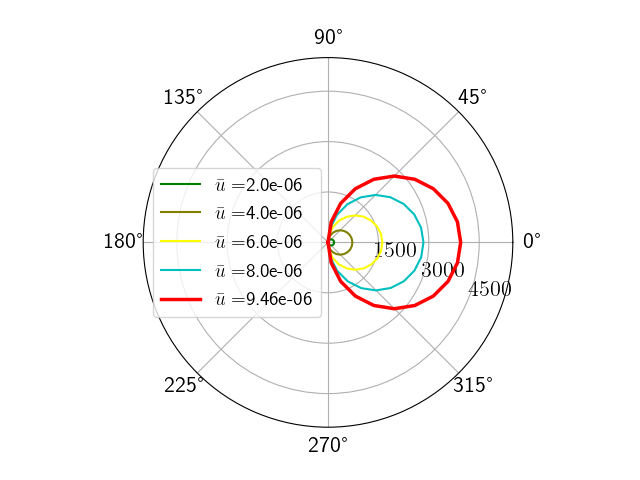}}
\quad
\subfigure[\texttt{no-tension}, $\alpha = 0$]{\includegraphics[width=0.35\linewidth]{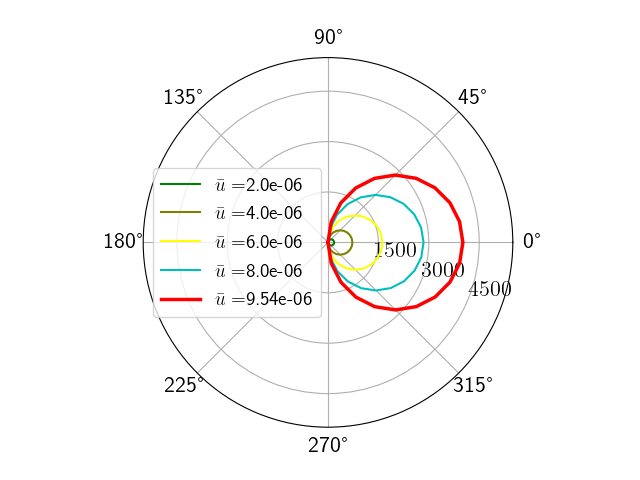}}
\vskip\baselineskip
\subfigure[\texttt{volumetric-deviatoric}, $\alpha = \pi/4$]{\includegraphics[width=0.35\linewidth]{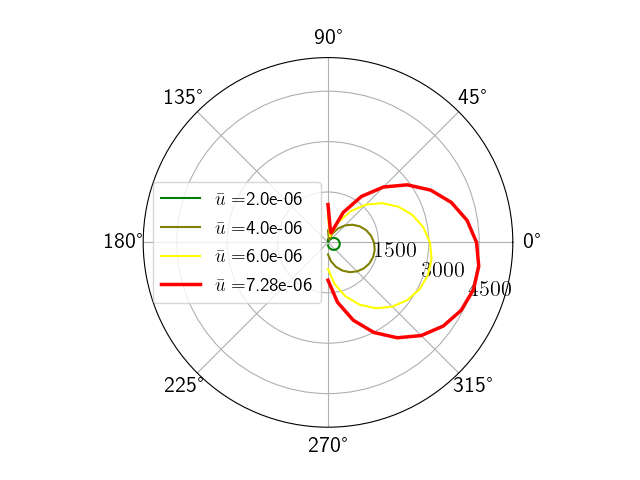}}
\quad
\subfigure[\texttt{no-tension}, $\alpha = \pi/4$]{\includegraphics[width=0.35\linewidth]{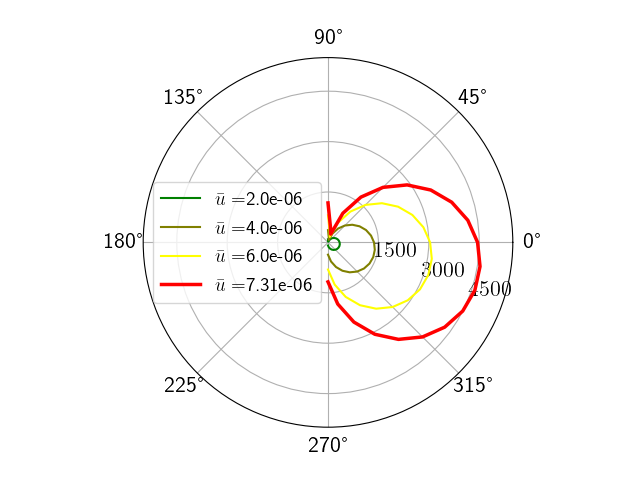}}
\vskip\baselineskip
\subfigure[\texttt{volumetric-deviatoric}, $\alpha = \pi/2$]{\includegraphics[width=0.35\linewidth]{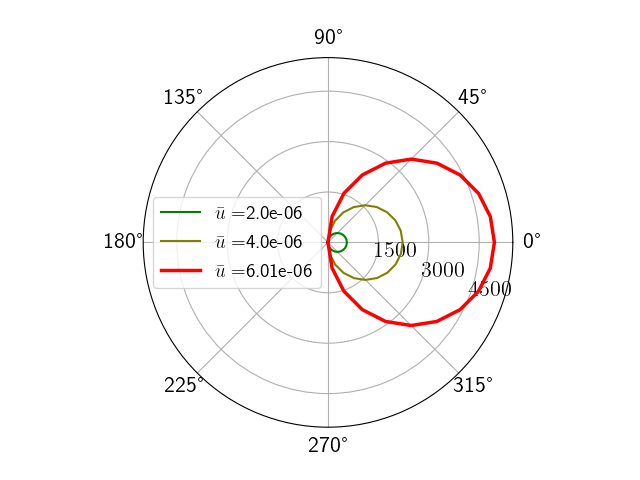}}
\quad
\subfigure[\texttt{no-tension}, $\alpha = \pi/2$]{\includegraphics[width=0.35\linewidth]{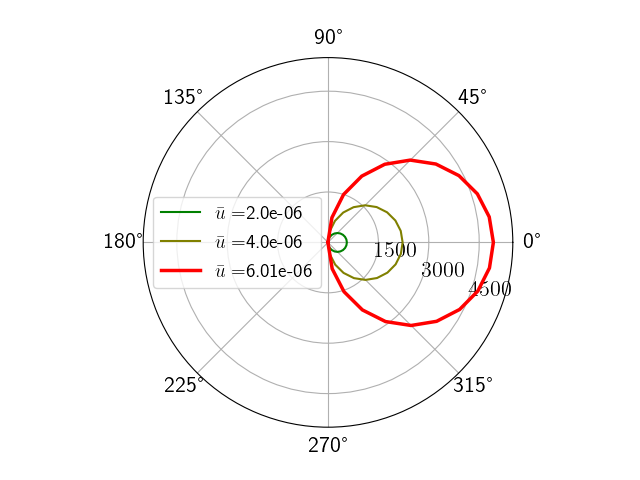}}
\caption{Under tensile loading for orthotropic {volumetric-deviatoric} and {no-tension}, we computed the energy release rate using the $G_{\theta}$ method with our proposed modifications. Several curves are plotted for a set of $\Bar{\bm{u}}$ where the red ones correspond with the stage of crack initiation. The peak of each curve approximately aligns with the direction of crack propagation, see Figure \ref{fig:tension-pf}. For $\alpha=0$ and $\alpha=\pi/2$, the peaks are along the horizontal line, whereas for $\alpha=\pm\pi/4$ the peaks deviate slightly from the $x$-axis. In both models, crack initiation occurs almost at the same loading for $\alpha=\pm\pi/4$ but in the opposite vertical directions, while the crack starts to propagate earlier for $\alpha=\pi/2$ compared to $\alpha=0$. This is expected due to the effect of orthotropy as the material is stiffer along the $y$-direction for $\alpha=\pi/2$.}
\label{fig:tension-gt-curve}
\end{figure}

A pre-notched square plate loaded in tension was computed. For orthotropic {volumetric-deviatoric} and {no-tension}, the fracture response at the last step $\Bar{\bm{u}}=2\times 10^{-5}$\,m is shown in Figure \ref{fig:tension-pf} with four material orientation angles $\alpha=0$, $\pm\pi/4$, and $\pi/2$. It was observed that both models generate similar crack paths. When the material principal axes coincide with the symmetric axes ($\alpha=0$ or $\alpha=\pi/2$), the crack moves straightforward. In other cases  ($\alpha=\pm\pi/4$), the crack slightly deviates from the straight path impacted by the orthotropic deformation. \newline
Figure~\ref{fig:tension-eps-ovd} depicts the magnitude of the strain tensor at the stage of the crack initiation for orthotropic {volumetric-deviatoric} with four material orientation angles. With $\alpha=0$ or $\alpha=\pi/2$, the deformation is symmetry around the $x$-axis whereas a non-symmetric behavior can be seen for $\alpha=\pm\pi/4$ where the material principal axes deviate from the symmetric axes. \newline
Figure \ref{fig:tension-gt-curve} shows the energy release rate computation using the $G_{\theta}$ method with our proposed modifications (see~\ref{App:G-theta})  for both models with four material orientation angles. Crack propagation should follow the direction of the maximum energy dissipation~\cite{Hakim2009,Chambolle2009kink}. And it is observed that for all cases the maximum energy release rate directions are approximately along the crack paths.
 Moreover, both split models have the identical crack propagation paths under the tensile loading because the strain energies are all tensile and crack-driving. Also observed is the crack propagates at the same loading with $\alpha=\pm\pi/4$, but earlier with $\alpha=\pi/2$ than $\alpha=0$. This is because different stiffness in the principal loading direction require different displacement loads to build up the necessary strain energy for propagating the crack. \newline

\subsection{Shear test}
\begin{figure}[htbp]
\centering %
\subfigure[\texttt{volumetric-deviatoric}, $\alpha = -\pi/4$]{\includegraphics[width=0.25\linewidth]{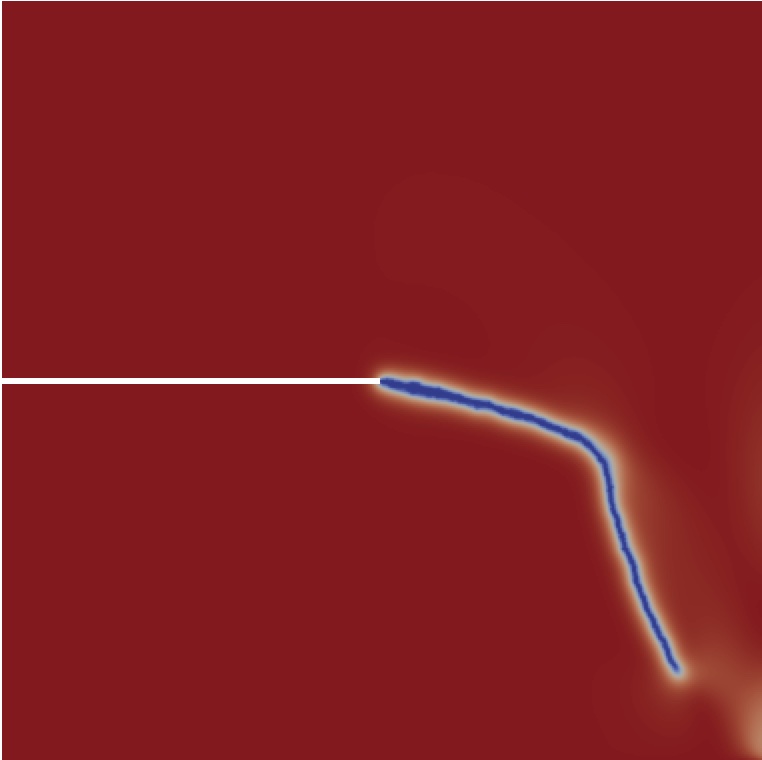}}
\quad
\subfigure[\texttt{no-tension}, $\alpha = -\pi/4$]{\includegraphics[width=0.25\linewidth]{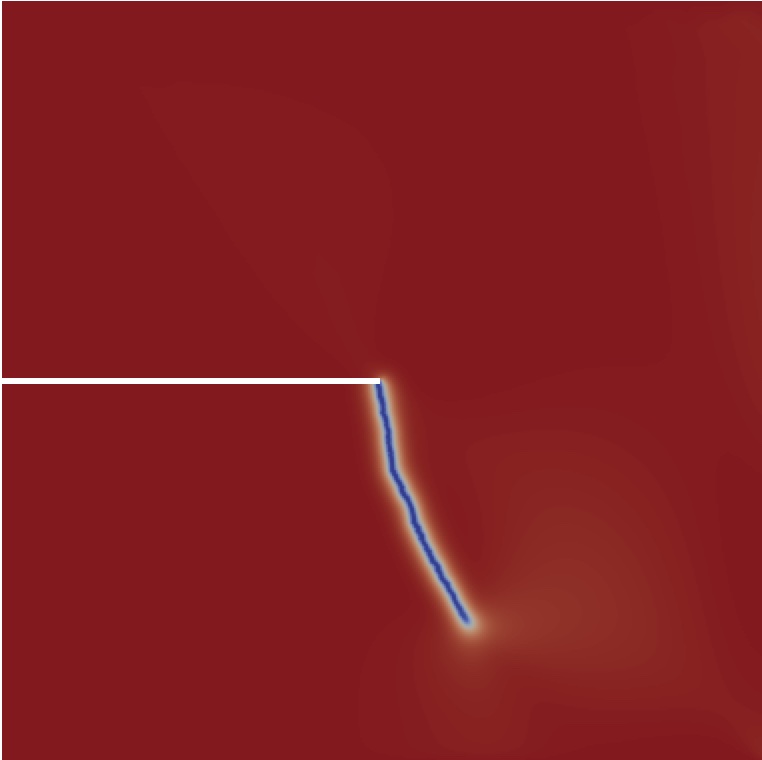}}
\vskip\baselineskip
\subfigure[\texttt{volumetric-deviatoric}, $\alpha = 0$]{\includegraphics[width=0.25\linewidth]{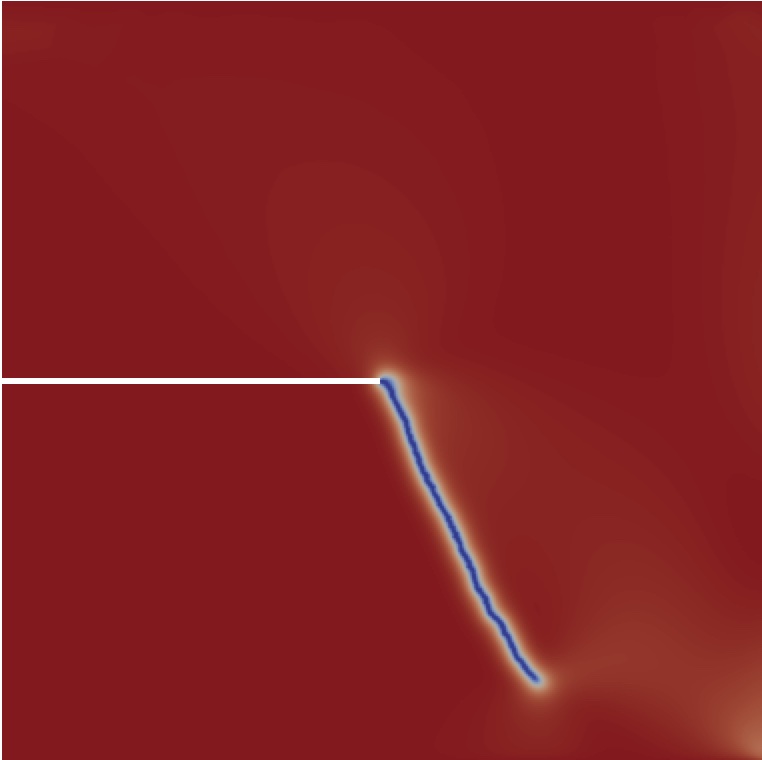}}
\quad
\subfigure[\texttt{no-tension}, $\alpha = 0$]{\includegraphics[width=0.25\linewidth]{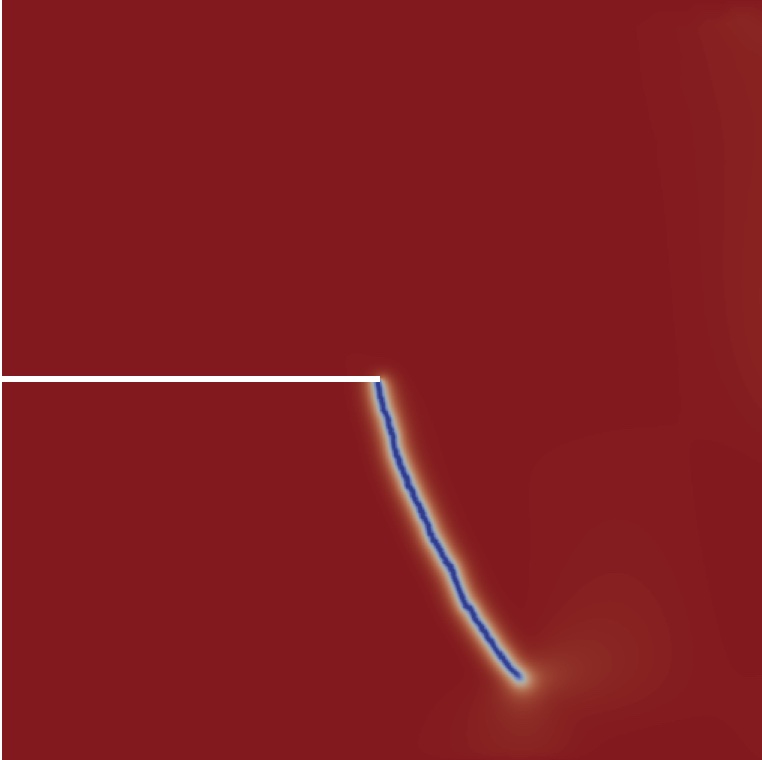}}
\vskip\baselineskip
\subfigure[\texttt{volumetric-deviatoric}, $\alpha = \pi/4$]{\includegraphics[width=0.25\linewidth]{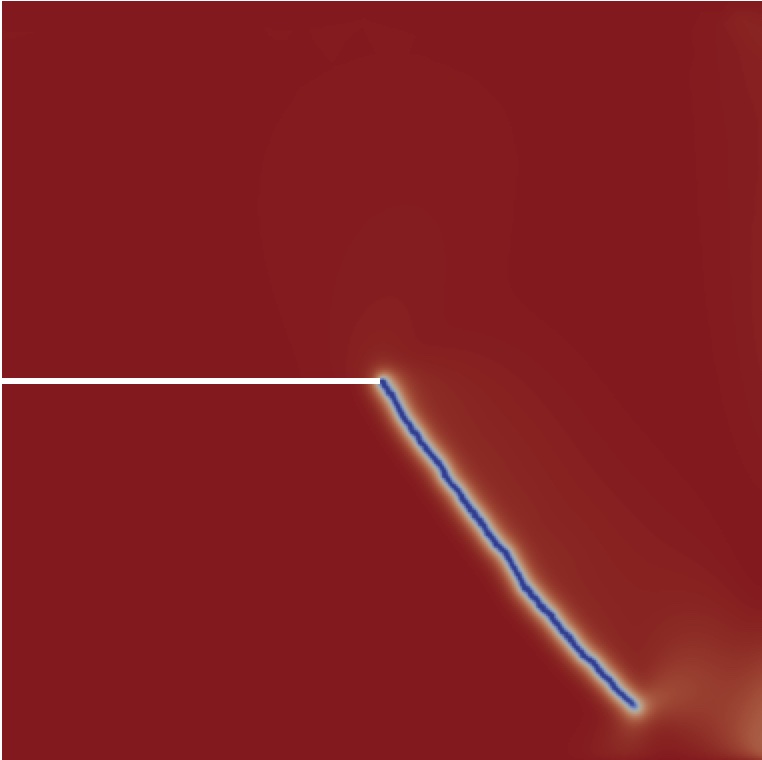}}
\quad
\subfigure[\texttt{no-tension}, $\alpha = \pi/4$]{\includegraphics[width=0.25\linewidth]{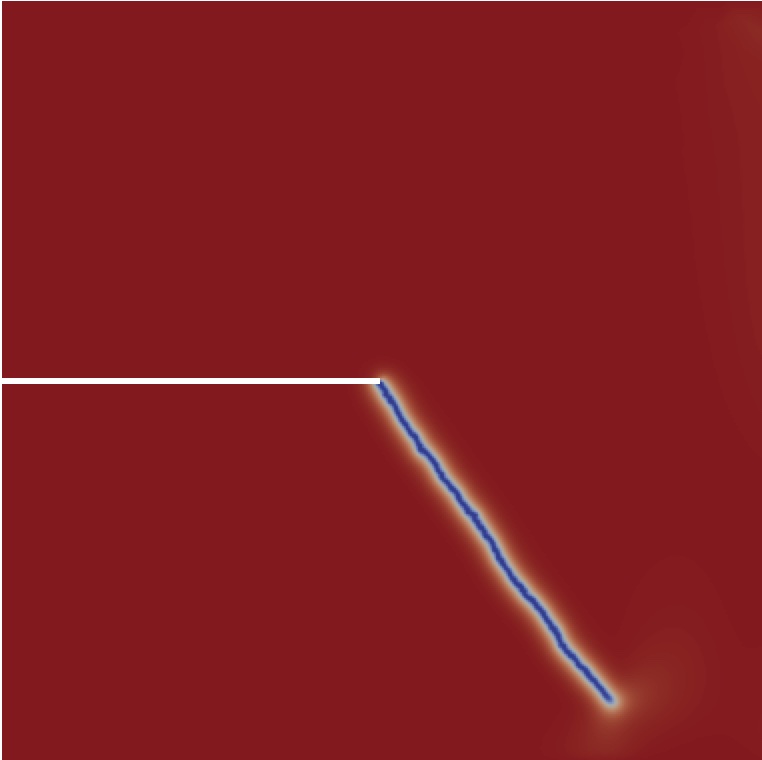}}
\vskip\baselineskip
\subfigure[\texttt{volumetric-deviatoric}, $\alpha = \pi/2$]{\includegraphics[width=0.25\linewidth]{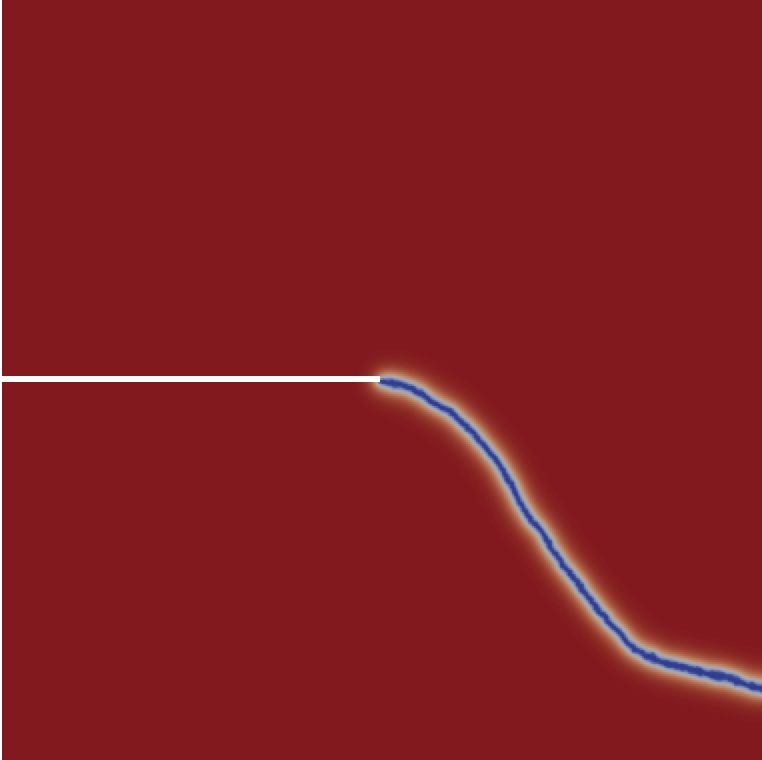}}
\quad
\subfigure[\texttt{no-tension}, $\alpha = \pi/2$]{\includegraphics[width=0.25\linewidth]{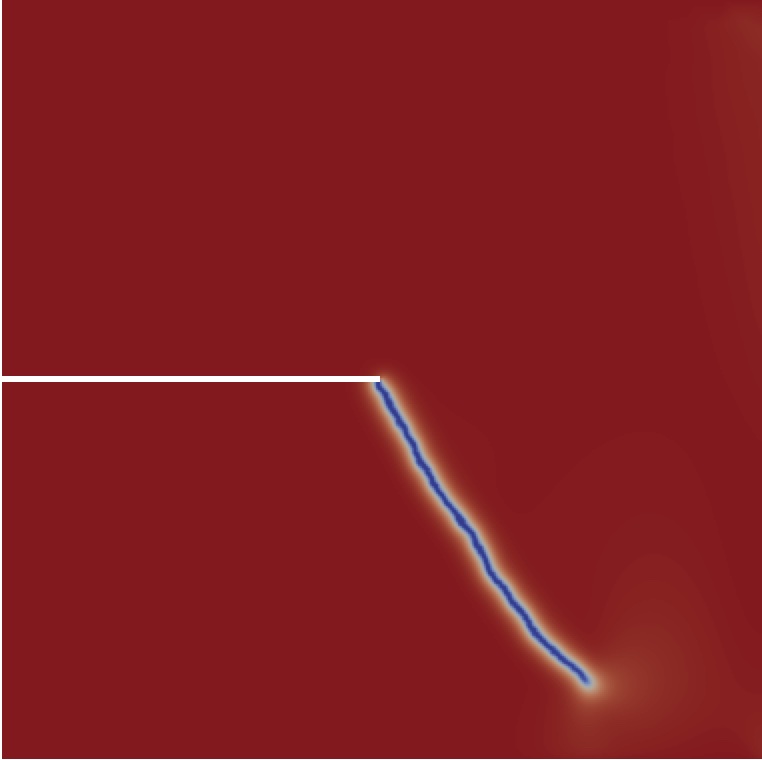}}
\caption{Cracked square plate under shear. Phase field contours of orthotropic {volumetric-deviatoric} and {no-tension} at the last step with $\Bar{\bm{u}} = 2\times 10^{-5}$ m. The same parameters were used for both models. Here, we observe that the crack paths are different depending on the model and the material orientation. For $\alpha=0$ and $\alpha=\pi/4$, however, both models generate relatively similar crack responses. 
With $\alpha=\pi/4$ there exist straight cracks that grow towards the lower right corner, almost along the material principal axis with lower stiffness. Also, with $\alpha=\pi/4$ the crack propagates earlier than other cases because of the effect of orthotropy. In this case, the material has the least stiffness along the preferred direction of crack growth under shear loading.}
\label{fig:shear-pf}
\end{figure}

\begin{figure}[htbp]
\centering %
\subfigure[\texttt{volumetric-deviatoric}, $\alpha = -\pi/4$]{\includegraphics[width=0.35\linewidth]{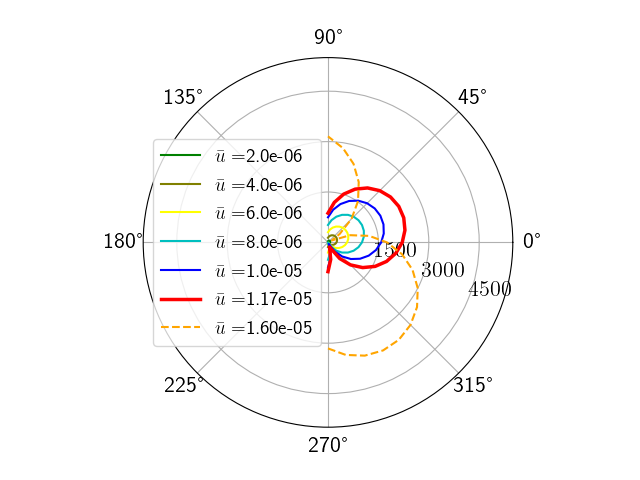}}
\quad
\subfigure[\texttt{no-tension}, $\alpha = -\pi/4$]{\includegraphics[width=0.35\linewidth]{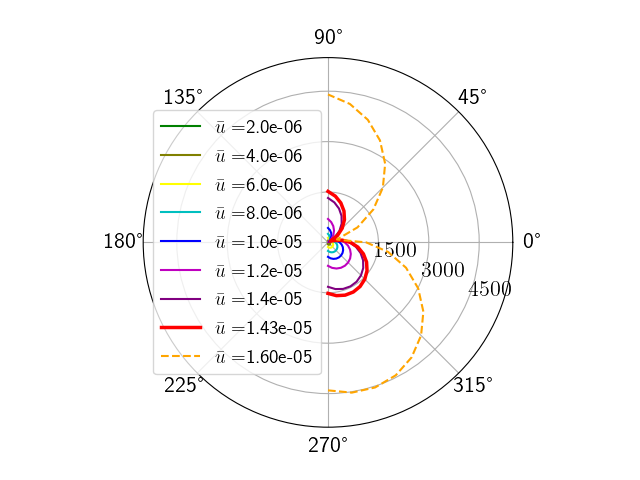}}
\vskip\baselineskip
\subfigure[\texttt{volumetric-deviatoric}, $\alpha = 0$]{\includegraphics[width=0.35\linewidth]{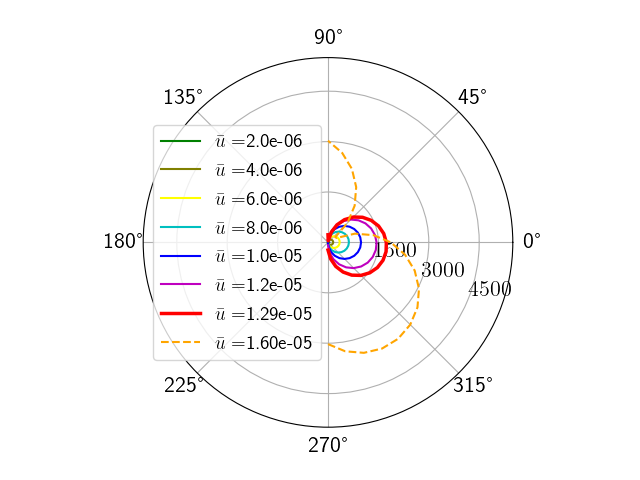}}
\quad
\subfigure[\texttt{no-tension}, $\alpha = 0$]{\includegraphics[width=0.35\linewidth]{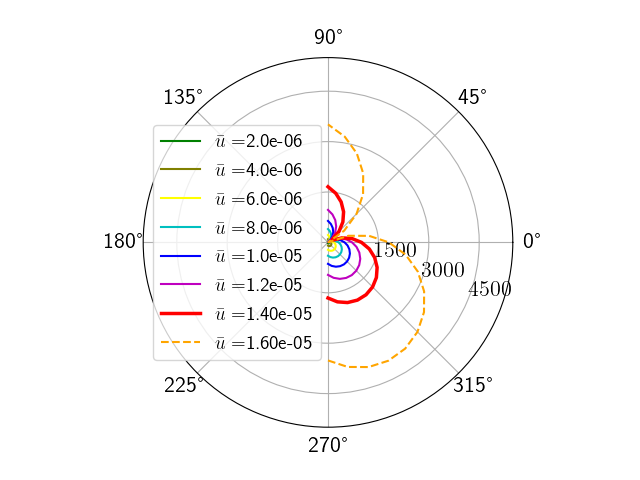}}
\vskip\baselineskip
\subfigure[\texttt{volumetric-deviatoric}, $\alpha = \pi/4$]{\includegraphics[width=0.35\linewidth]{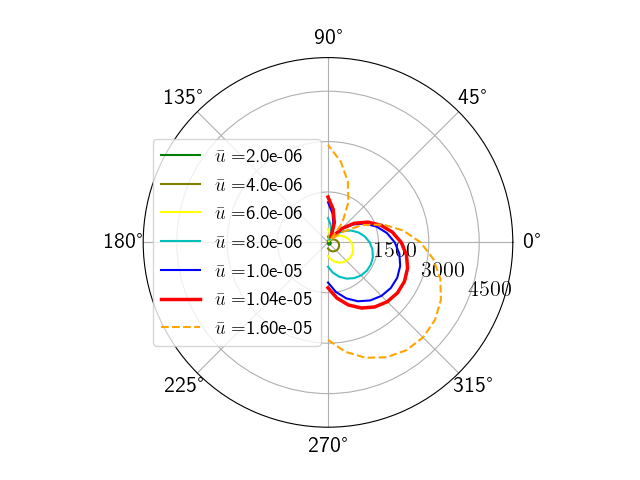}}
\quad
\subfigure[\texttt{no-tension}, $\alpha = \pi/4$]{\includegraphics[width=0.35\linewidth]{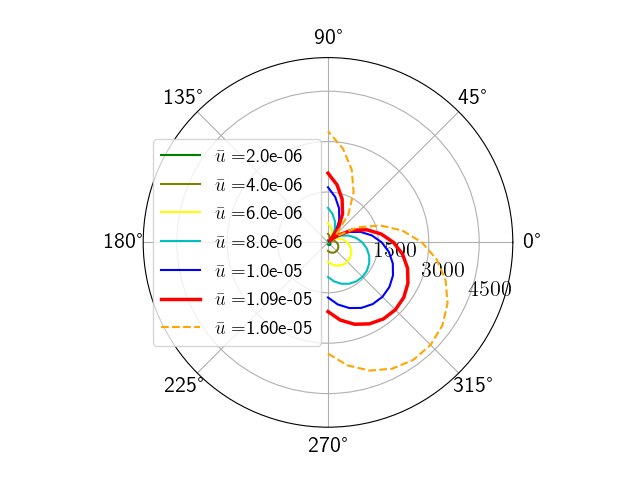}}
\vskip\baselineskip
\subfigure[\texttt{volumetric-deviatoric}, $\alpha = \pi/2$]{\includegraphics[width=0.35\linewidth]{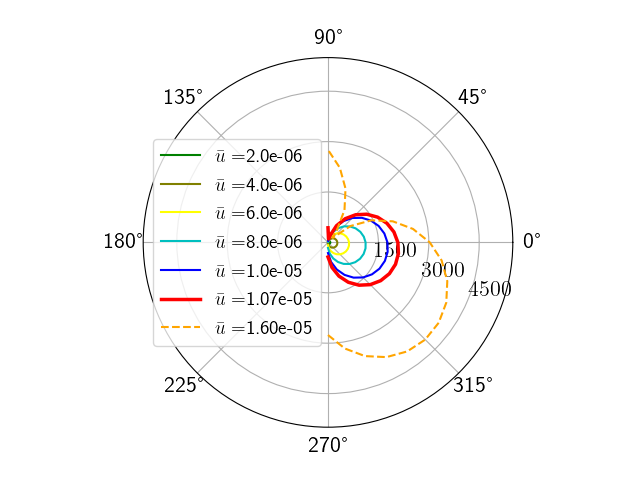}}
\quad
\subfigure[\texttt{no-tension}, $\alpha = \pi/2$]{\includegraphics[width=0.35\linewidth]{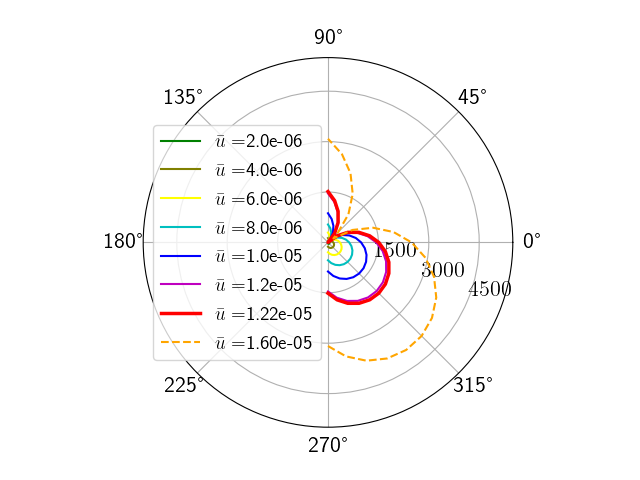}}
\caption{Under shear loading for orthotropic {volumetric-deviatoric} and {no-tension}, we computed the energy release rate using the $G_{\theta}$ method with our proposed modifications. Several curves are plotted for a set of $\Bar{\bm{u}}$ where the red curves correspond with the onset of crack propagation, and the orange ones refer to a fixed step thereafter for all samples, $\Bar{\bm{u}}=1.6\times 10^{-5}$ m,  where the crack is propagating towards the bottom-right corner. Like the tensile example, for most cases, the maximum energy release rate of each case approximately aligns with the direction of crack propagation (Figure \ref{fig:shear-pf}). An exception is the case of $\alpha=-\pi/4$ for {volumetric-deviatoric} where the peak of the red curve deviates from the direction of crack propagation, whereas the orange curve fits better. Also observed is that crack initiation occurs earlier in {no-tension} for all cases.}
\label{fig:shear-gt-curve}
\end{figure}

\begin{figure}[htbp]
\centering %
\subfigure[\texttt{volumetric-deviatoric}, tension]{\includegraphics[width=0.35\linewidth]{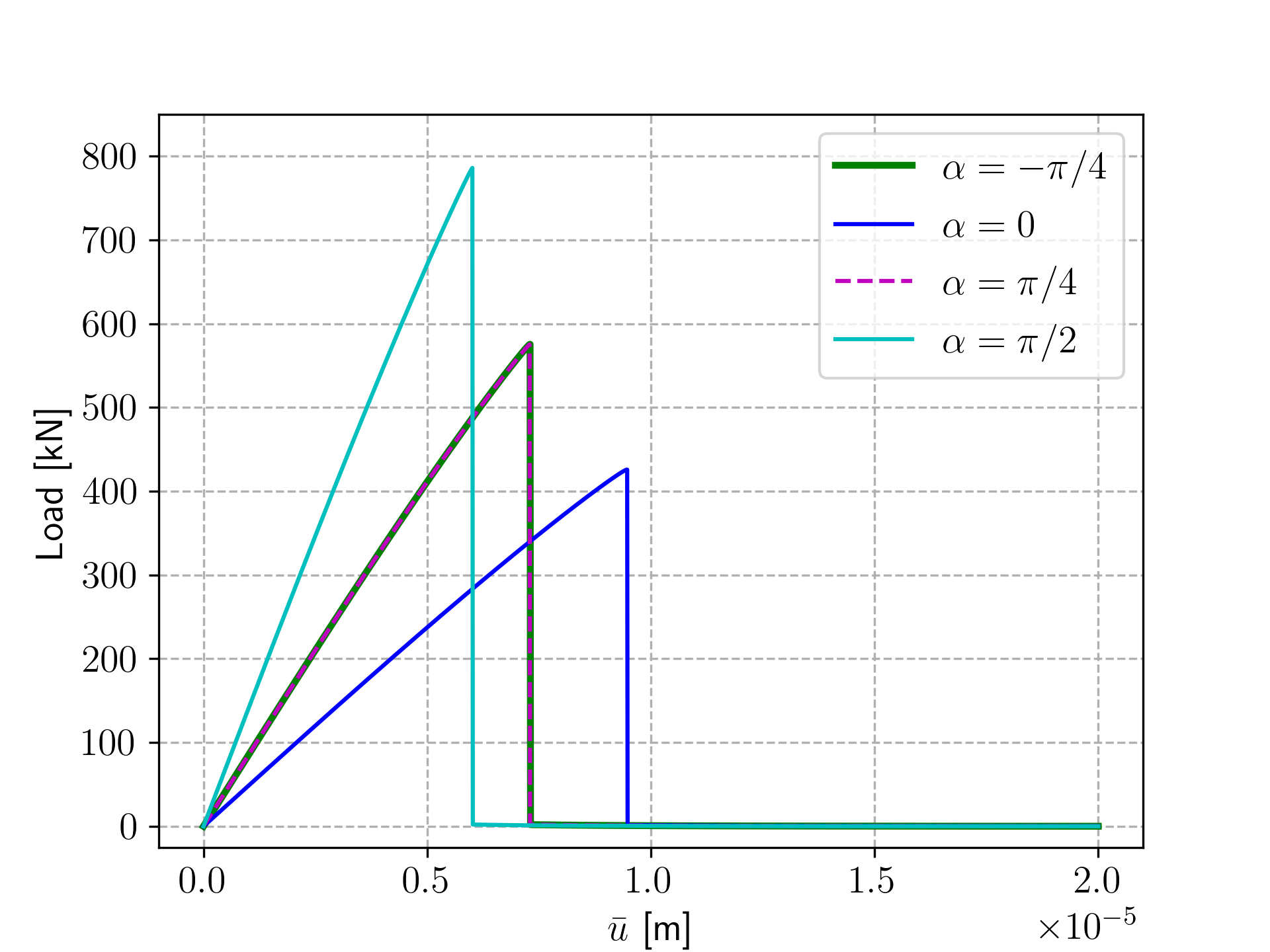}}
\quad
\subfigure[\texttt{no-tension}, tension]{\includegraphics[width=0.35\linewidth]{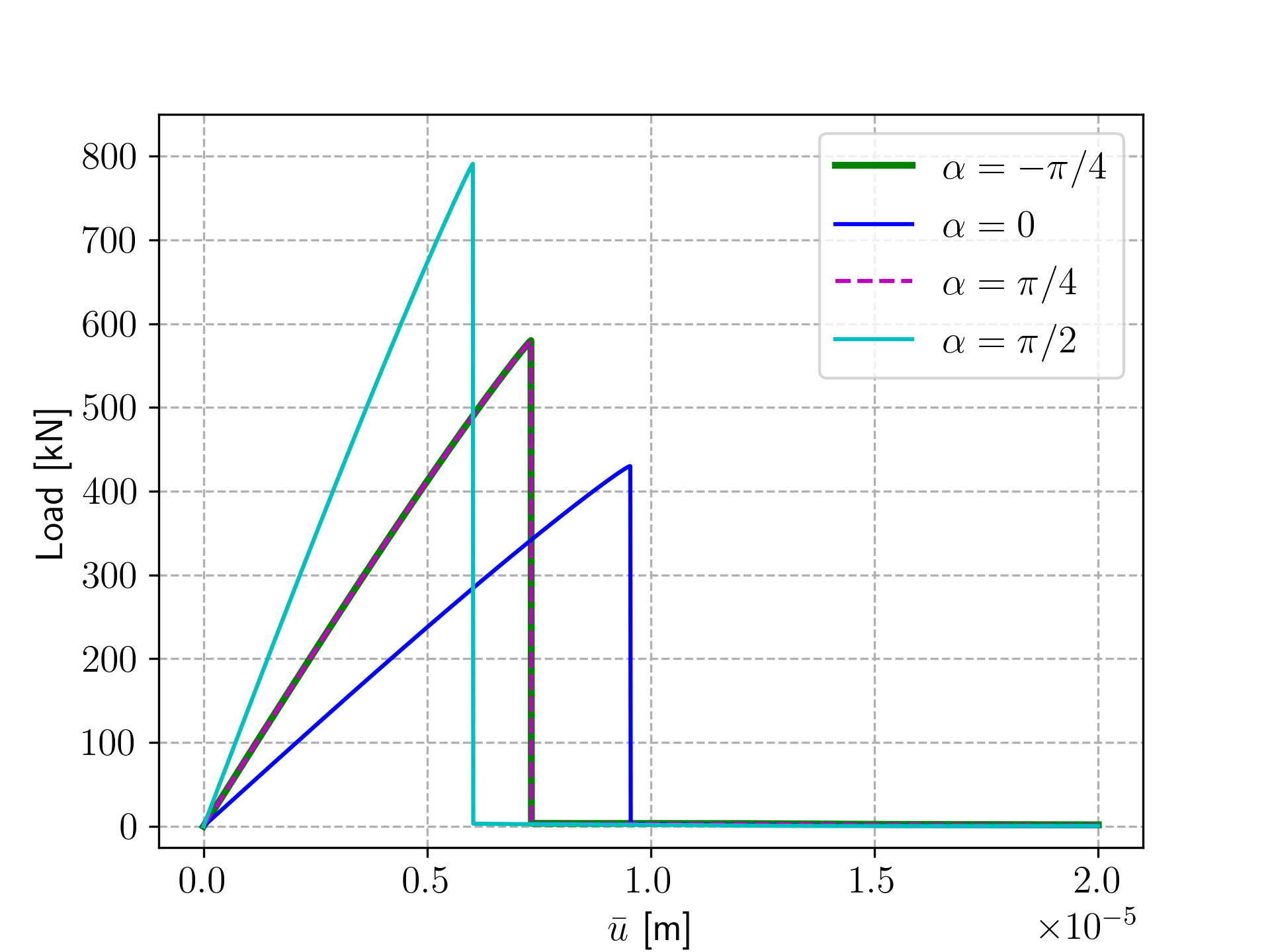}}
\vskip\baselineskip
\subfigure[\texttt{volumetric-deviatoric}, shear]{\includegraphics[width=0.35\linewidth]{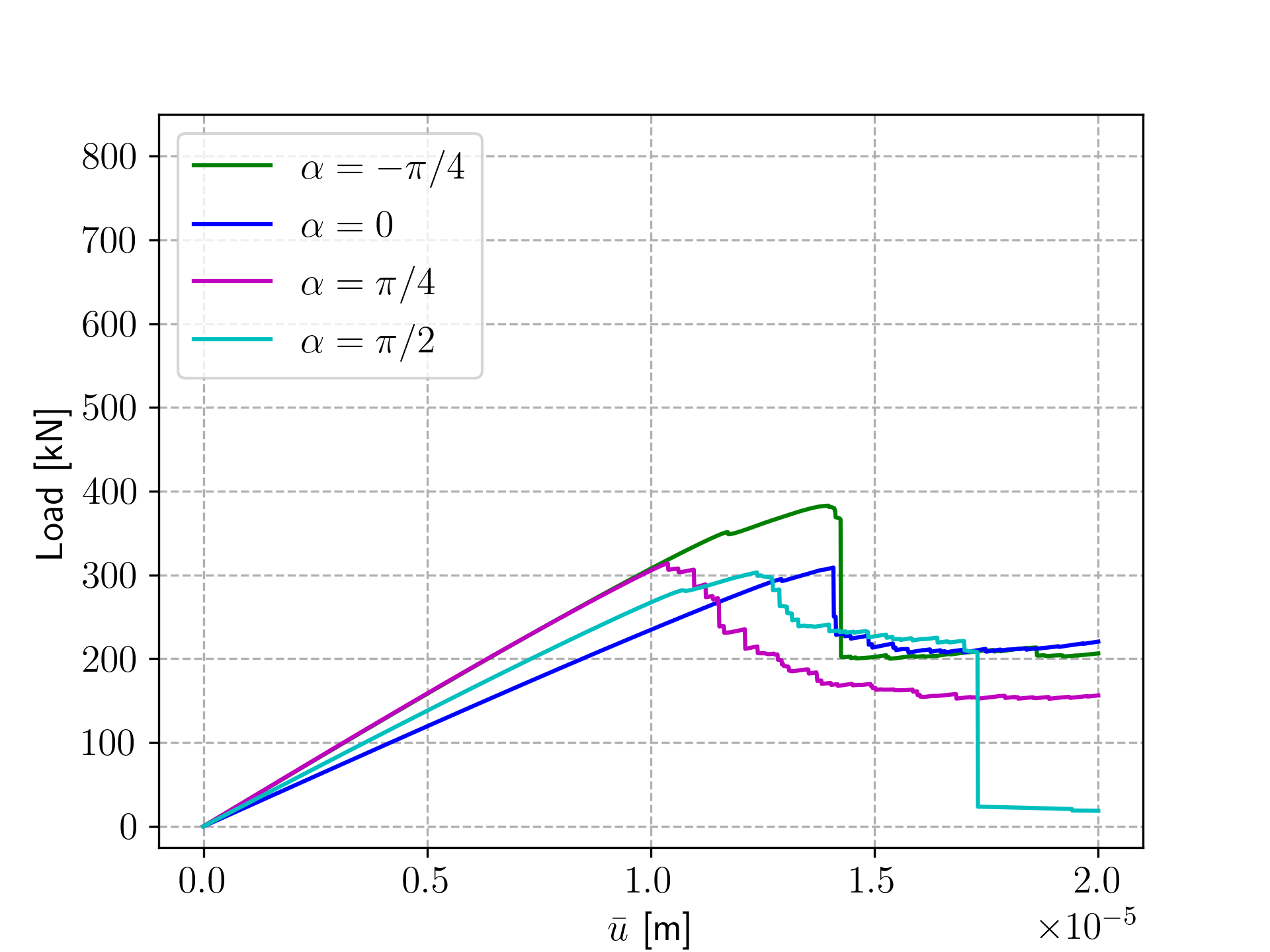}}
\quad
\subfigure[\texttt{no-tension}, shear]{\includegraphics[width=0.35\linewidth]{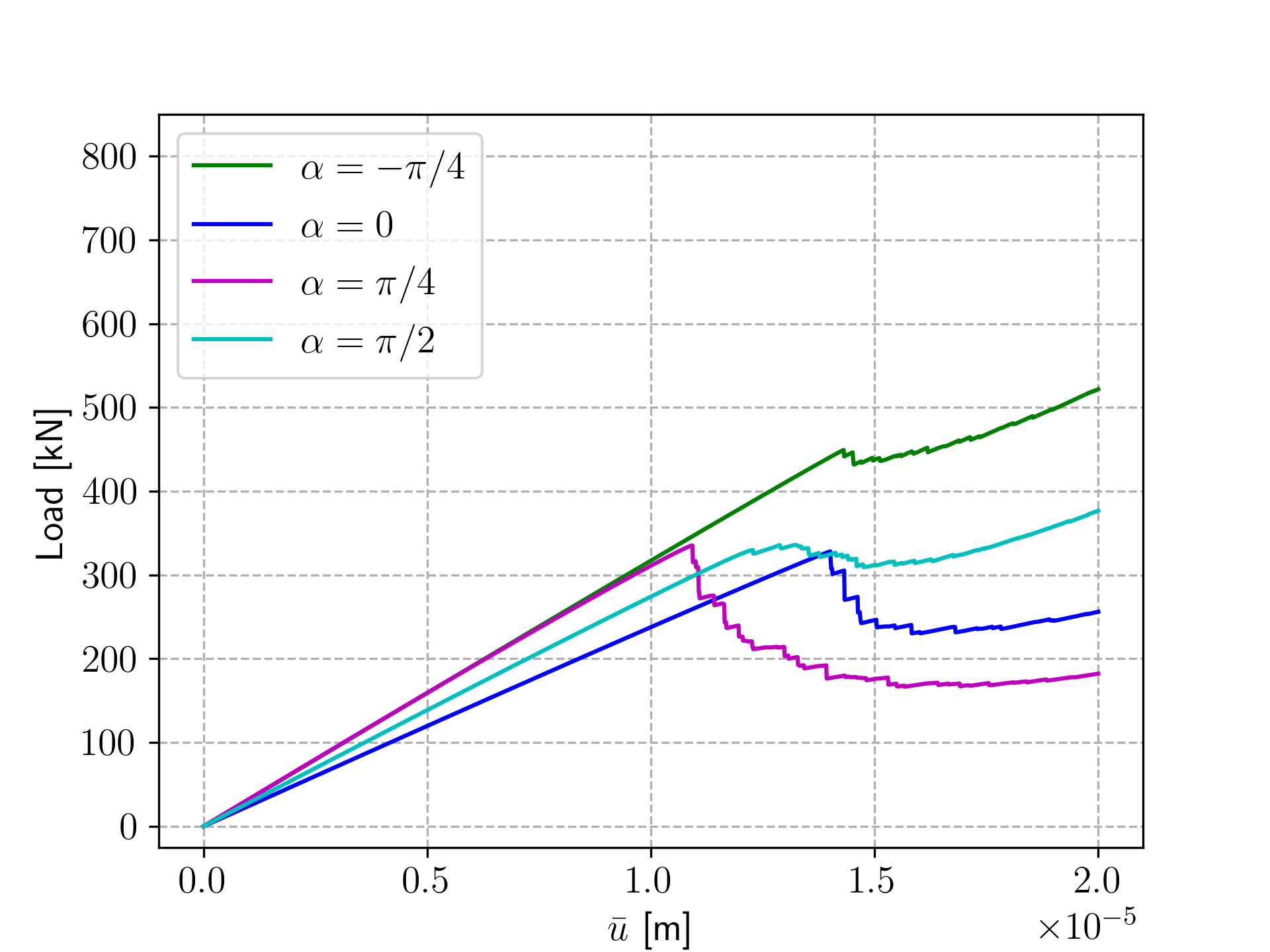}}
\caption{Load deflection curves for both models under (a)--(b) tensile and (c)--(d) shear loading. For (a)--(b) a sudden loss of the stiffness indicates that the crack grows in an unstable manner. These load-displacment results show that the crack propagates at the same loading for $\alpha=\pm\pi/4$. For $\alpha=\pi/2$ the crack propagates earlier compared to $\alpha=0$ because the crack propagation direction aligns with the weaker material direction. For (c), the global stiffness is ``lost" for the case of $\alpha=\pi/2$ in which the crack reaches to the edge. Also observed for all cases in (c)--(d) is that there is a drop in the load at the onset of crack propagation.}
\label{fig:load-disp}
\end{figure}

Figure \ref{fig:shear-pf} depicts the crack paths for orthotropic {volumetric-deviatoric} and {no-tension} under shear loading with four material orientation angles ($\alpha=-\pi/4$, 0, $\pi/4$, and $\pi/2$). In contrary to the tensile test, herein a significant difference in crack response is observed between the two models for all cases.
Figure \ref{fig:shear-gt-curve} plots the maximum energy release rate for both models with four material orientation angles. In most cases, the peaks are approximately along the crack propagation paths. Note, however, that the maximum energy release rate computed is higher than the critical value which may require further investigation.
Figures \ref{fig:load-disp} show the load-deflection curves for the tensile examples (a--b) and the shear examples (c--d). All tensile examples result in a global stiffness loss, as the crack propagates through the sample. For the shear examples, only in one case ({volumetric-deviatoric} with $\alpha=\pi/2$), the crack reaches to the edge. In other cases, as soon as the crack propagates, the load drops. However this drop is not high and is recovered shortly in some cases.)

    \pagebreak
    \section{Conclusion}
    In this paper, we have proposed a method to account for orthotropy/anisotropy of materials when modeling tension-compression asymmetry in their crack response. Two decomposition models, {volumetric-deviatoric} and {no-tension}, were extended to capture arbitrary types of anisotropy in the constitutive behavior in a three-dimensional setting. 
A major characteristic of the formulation is that it keeps the variational formalism in partitioning the constitutive model, honoring the orthogonality condition. \newline
The two models presented in this work were evaluated and compared numerically in a square plate with a pre-notched crack loaded in tension and shear. 
Under the tensile loading, the result indicates that both models have approximately similar responses. Testing different material orientation angles, the effect of orthotropy on the material response is clearly observed. To better interpret the results, we investigated the energy release rate using the $G_\theta$ method with our proposed modifications.
And a significant difference in the resulting crack paths iss observed between the two models under shear loading. However, at this point it is difficult to judge which model provides a more reliable outcome. \newline
The results presented herein highlight the importance of how to decompose the constitutive model of materials with anisotropic nature in the phase field approach to brittle fracture. 
As opposed to imposing anisotropic behaviors to the surface energy (fracture toughness), adapting the strain energy leads to `weak' anisotropy~\cite{Li2019aniso}.
However, some materials evidently exhibit anisotropic elastic deformation and this effect should not be ignored.
Future studies should include `strong' anisotropy by adding the directionality to the surface energy, and comparisons against experiments to determine whether both strong and weak anisotropy or how much of each is required.

    \section*{Acknowledgements}
    VZR, OK, and TN gratefully acknowledge the funding provided by the German Federal Ministry of Education and Research (BMBF) for the GeomInt2 project (grand number 03G0899D) as well as the support of the Project Management J{\"u}lich (PtJ). MM and KY's contributions are funded by the Helmholtz Association (grant number SO-093) through the iCROSS-Project, EURAD, the European Joint Programme on Radioactive Waste Management through the DONUT and MAGIC work packages (Grant Agreement No 847593), and the Deutsche Forschungsgemeinschaft (DFG, German Research Foundation) through the HIGHER project (grant number PA 3451/1‐1).
    The authors thank the Earth System Modelling Project (ESM) for funding this work by providing computing time on the ESM partition of the supercomputer JUWELS at the J{\"u}lich Supercomputing Centre (JSC). The authors are very grateful to the OpenGeoSys developer core team for technical advice in coding and software implementation support.
    \appendix
    \section{Projection tensors in two/three dimensions}
\label{App}
{For the two-dimensional case, the following closed-form expressions are adapted.}
\begin{subequations}\label{Eq:mas-eps-deriv-2d}
    \begin{align}
        \mathbb{D}^\text{t} := \frac{\partial\tilde{\bm{\varepsilon}}^\text{t}}{\partial\tilde{\bm{\varepsilon}}} = 
        \begin{cases}
            \dfrac{a_1 - a_2}{\tilde{\varepsilon}_1 - \tilde{\varepsilon}_2}\Bigl[\mathds{1}_S - \bm{M}_1 \otimes \bm{M}_1 - \bm{M}_2 \otimes \bm{M}_2 \Bigr] \\
            + H(\tilde{\varepsilon}_{1})\bm{M}_{1}\otimes\bm{M}_{1} + H(\tilde{\varepsilon}_{2})\bm{M}_{2}\otimes\bm{M}_{2},
            & \text{if}~ \tilde{\varepsilon}_1\neq \tilde{\varepsilon}_2 \\
            H(\tilde{\varepsilon}_1)\mathds{1}_S, 
            & \text{if}~ \tilde{\varepsilon}_1 = \tilde{\varepsilon}_2
        \end{cases} \\
        \mathbb{D}^\text{c} := \frac{\partial\tilde{\bm{\varepsilon}}^\text{c}}{\partial\tilde{\bm{\varepsilon}}} = 
        \begin{cases}
            \dfrac{b_1 - b_2}{\tilde{\varepsilon}_1 - \tilde{\varepsilon}_2}\Bigl[\mathds{1}_S - \bm{M}_1 \otimes \bm{M}_1 - \bm{M}_2 \otimes \bm{M}_2 \Bigr] \\
            + H(-\tilde{\varepsilon}_{1})\bm{M}_{1}\otimes\bm{M}_{1} + H(-\tilde{\varepsilon}_{2})\bm{M}_{2}\otimes\bm{M}_{2},
            & \text{if}~ \tilde{\varepsilon}_1\neq \tilde{\varepsilon}_2 \\
            H(-\tilde{\varepsilon}_1)\mathds{1}_S,
            & \text{if}~ \tilde{\varepsilon}_1 = \tilde{\varepsilon}_2
        \end{cases}
    \end{align}
\end{subequations}
where $\mathds{1}_S$ is the fourth-order tensor defined by:
\begin{equation*}
    \begin{aligned}
        (\mathds{1}_S)_{ijkl} = \frac12\bigl(\delta_{ik}\delta_{jl} + \delta_{il}\delta_{jk}\bigr).
    \end{aligned}
\end{equation*}
To derive the projection tensors ($\mathbb{D}^\text{t}$ and $\mathbb{D}^\text{c}$) with a three-dimensional basis in Kelvin-matrix notation, we proceed as follows. With  \eqref{Eq:proj-dir-dev}, let the original coordinate system coincide with the material principal axes $\bm{n}_i$s.
Thus, the only non-zero components of $\mathbb{D}^{t}$ and $\mathbb{D}^{c}$ on the basis $\bm{n}_i$ are given by:
\begin{equation}\label{Eq:nonzero-dir-dev}
    \begin{aligned}
        D^\text{t}_{iiii} = H(\tilde{\varepsilon}_i),
        \quad D^\text{c}_{iiii} = H(-\tilde{\varepsilon}_i),
    \end{aligned}
\end{equation}
\begin{equation*}
    \begin{aligned}
        D^\text{t}_{ijij} = 
        \begin{cases}
            \frac12 H(\tilde{\varepsilon}_i) &
            \tilde{\varepsilon}_i=\tilde{\varepsilon}_j,\quad i\neq j, \\
            \frac12 \bigl\{H(\tilde{\varepsilon}_j)\tilde{\varepsilon}_j - H(\tilde{\varepsilon}_i)\tilde{\varepsilon}_i \bigr\}/(\tilde{\varepsilon}_j - \tilde{\varepsilon}_i) & \tilde{\varepsilon}_i\neq \tilde{\varepsilon}_j,\quad i\neq j,
        \end{cases} \\
        D^\text{c}_{ijij} = 
        \begin{cases}
            \frac12 H(-\tilde{\varepsilon}_i) &
            \tilde{\varepsilon}_i = \tilde{\varepsilon}_j,\quad i\neq j, \\
            \frac12\bigl\{H(-\tilde{\varepsilon}_j)\tilde{\varepsilon}_j - H(-\tilde{\varepsilon}_i)\tilde{\varepsilon}_i \bigr\}/(\tilde{\varepsilon}_j - \tilde{\varepsilon}_i) & \tilde{\varepsilon}_i\neq \tilde{\varepsilon}_j,\quad i\neq j.
        \end{cases}
    \end{aligned}
\end{equation*}
The following projection tensors are then  obtained for different cases (in Kelvin-matrix notation):
\paragraph{Case \RNum{1}}
$\tilde{\bm{\varepsilon}}^\text{t}=(\tilde{\varepsilon}_1,\tilde{\varepsilon}_2,\tilde{\varepsilon}_3)$ and $\tilde{\bm{\varepsilon}}^\text{c}=(0,0,0)$:
\begin{equation*}
    \begin{aligned}
        \undertilde{\bm{D}}^\text{t} = 
        \begin{bmatrix}
            1 & 0 & 0 & 0 & 0 & 0 \\
            & 1 & 0 & 0 & 0 & 0 \\
            &  & 1 & 0 & 0 & 0 \\
            &  &  & 1 & 0 & 0 \\
            &  &  &  & 1 & 0\\
            &  &  &  &  & 1
        \end{bmatrix}, \quad
        \undertilde{\bm{D}}^\text{c} = 
        \begin{bmatrix}
            0 & 0 & 0 & 0 & 0 & 0 \\
            & 0 & 0 & 0 & 0 & 0 \\
            &  & 0 & 0 & 0 & 0 \\
            &  &  & 0 & 0 & 0 \\
            &  &  &  & 0 & 0 \\
            &  &  &  &  & 0
        \end{bmatrix}.
    \end{aligned}
\end{equation*}
\paragraph{Case \RNum{2}} $\tilde{\bm{\varepsilon}}^\text{t}=(\tilde{\varepsilon}_1,\tilde{\varepsilon}_2,0)$ and $\tilde{\bm{\varepsilon}}^\text{c}=(0,0,\tilde{\varepsilon}_3)$:
\begin{equation*}
    \begin{aligned}
        \undertilde{\bm{D}}^\text{t} = 
        \begin{bmatrix}
            1 & 0 & 0 & 0 & 0 & 0 \\
            & 1 & 0 & 0 & 0 & 0 \\
            &  & 0 & 0 & 0 & 0 \\
            &  &  & 1 & 0 & 0 \\
            &  &  &  & \tilde{\varepsilon}_1/(\tilde{\varepsilon}_1-\tilde{\varepsilon}_3) & 0 \\
            &  &  &  &  & \tilde{\varepsilon}_2/(\tilde{\varepsilon}_2-\tilde{\varepsilon}_3)
        \end{bmatrix}, \quad \undertilde{\bm{D}}^\text{c} = 
        \begin{bmatrix}
            0 & 0 & 0 & 0 & 0 & 0 \\
            & 0 & 0 & 0 & 0 & 0 \\
            &  & 1 & 0 & 0 & 0 \\
            &  &  & 0 & 0 & 0 \\
            &  &  &  & -\tilde{\varepsilon}_3/(\tilde{\varepsilon}_1-\tilde{\varepsilon}_3) & 0 \\
            &  &  &  &  & -\tilde{\varepsilon}_3/(\tilde{\varepsilon}_2-\tilde{\varepsilon}_3)
        \end{bmatrix}.
    \end{aligned}
\end{equation*}
\paragraph{Case \RNum{3}}
$\tilde{\bm{\varepsilon}}^\text{t}=(\tilde{\varepsilon}_1,0,0)$ and $\tilde{\bm{\varepsilon}}^\text{c}=(0,\tilde{\varepsilon}_2,\tilde{\varepsilon}_3)$:
\begin{equation*}
    \begin{aligned}
        \undertilde{\bm{D}}^\text{t} = 
        \begin{bmatrix}
            1 & 0 & 0 & 0 & 0 & 0 \\
            & 0 & 0 & 0 & 0 & 0 \\
            &  & 0 & 0 & 0 & 0 \\
            &  &  & \tilde{\varepsilon}_1/(\tilde{\varepsilon}_1-\tilde{\varepsilon}_2) & 0 & 0 \\
            &  &  &  & \tilde{\varepsilon}_1/(\tilde{\varepsilon}_1-\tilde{\varepsilon}_3) & 0 \\
            &  &  &  &  & 0
        \end{bmatrix}, \quad \undertilde{\bm{D}}^\text{c} = 
        \begin{bmatrix}
            0 & 0 & 0 & 0 & 0 & 0 \\
            & 1 & 0 & 0 & 0 & 0 \\
            &  & 1 & 0 & 0 & 0 \\
            &  &  & -\tilde{\varepsilon}_2/(\tilde{\varepsilon}_1-\tilde{\varepsilon}_2) & 0 & 0 \\
            &  &  &  & -\tilde{\varepsilon}_3/(\tilde{\varepsilon}_1-\tilde{\varepsilon}_3) & 0\\
            &  &  &  &  & 1
        \end{bmatrix}.
    \end{aligned}
\end{equation*}
\paragraph{Case \RNum{4}}
$\tilde{\bm{\varepsilon}}^\text{t}=(0,0,0)$ and $\tilde{\bm{\varepsilon}}^\text{c}=(\tilde{\varepsilon}_1,\tilde{\varepsilon}_2,\tilde{\varepsilon}_3)$:
\begin{equation*}
    \begin{aligned}
        \undertilde{\bm{D}}^\text{t} = 
        \begin{bmatrix}
            0 & 0 & 0 & 0 & 0 & 0 \\
            & 0 & 0 & 0 & 0 & 0 \\
            &  & 0 & 0 & 0 & 0 \\
            &  &  & 0 & 0 & 0 \\
            &  &  &  & 0 & 0 \\
            &  &  &  &  & 0
        \end{bmatrix}, \quad \undertilde{\bm{D}}^\text{c} = 
        \begin{bmatrix}
            1 & 0 & 0 & 0 & 0 & 0 \\
            & 1 & 0 & 0 & 0 & 0 \\
            &  & 1 & 0 & 0 & 0 \\
            &  &  & 1 & 0 & 0 \\
            &  &  &  & 1 & 0\\
            &  &  &  &  & 1
        \end{bmatrix}.
    \end{aligned}
\end{equation*}
The projection tensors in Kelvin-matrix notation, $\undertilde{\bm{D}}^\text{t}$ and $\undertilde{\bm{D}}^\text{c}$ can also be derived for an arbitrary orthonormal basis with a transformation as follows
\begin{equation}\label{Eq:project-tensor-transform}
    \begin{aligned}
        \undertilde{\bm{D}}^{\prime \text{t}} = \bm{P}\undertilde{\bm{D}}^{\text{t}}\bm{P}^T, \quad \undertilde{\bm{D}}^{\prime \text{c}} = \bm{P}\undertilde{\bm{D}}^{\text{c}}\bm{P}^T,
    \end{aligned}
\end{equation}
where $\undertilde{\bm{D}}^{\prime \text{t}}$ and $\undertilde{\bm{D}}^{\prime \text{c}}$ are represented in the new basis.

    \section{Energy release rate computation: $G_\theta$ method}\label{App:G-theta}

\begin{figure}[htbp]
\centering %
\subfigure[]{\includegraphics[width=0.35\linewidth]{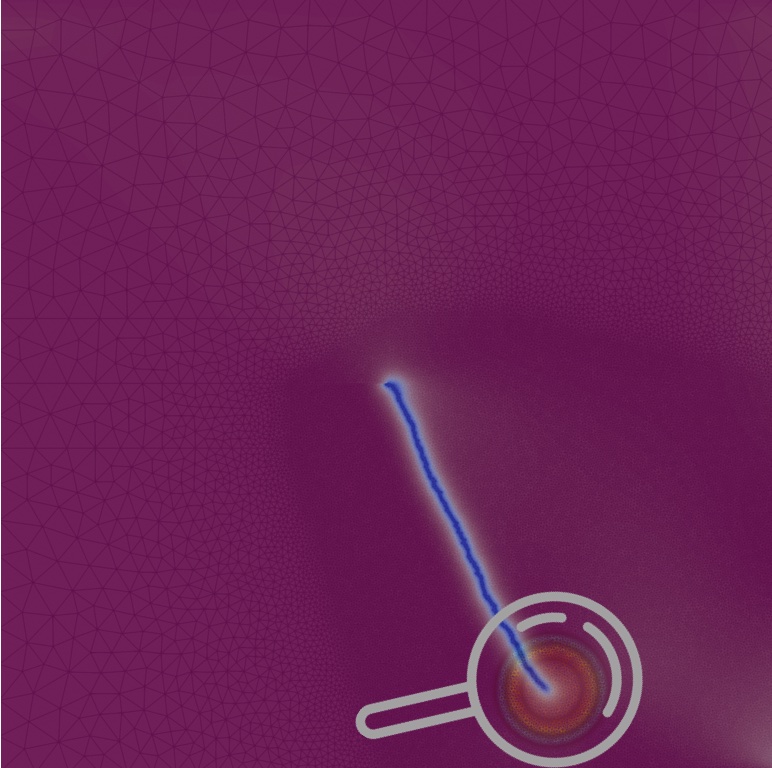}}
\quad
\subfigure[]{\includegraphics[width=0.25\linewidth]{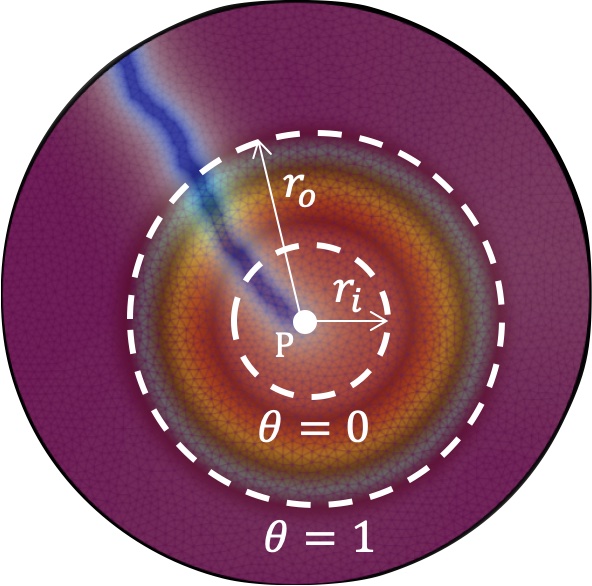}}\\
\subfigure{\includegraphics[width=0.3\linewidth]{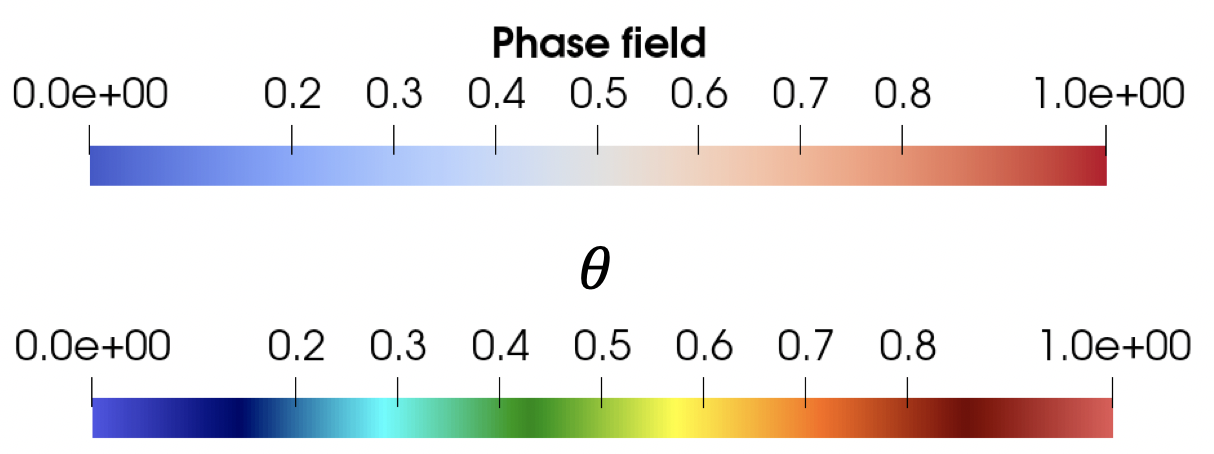}}

\caption{Phase field and $\theta$ profile  for \texttt{volumetric-deviatoric},  $\alpha=0$ at $t_f$ under shear loading.  We use virtual perturbation of $\theta$ to compute energy release rate using $G_{\theta}$. The  $\theta$ value is 1 inside of $B_{r_{i}}(P)$,  0 outside, and a linear interpolation in between. We set $r_{i}=4\ell$ and $r_{o}=2.5r_{i}$.}
\label{fig:Gtheta_schematic}
\end{figure}

The $G_{\theta}$ method is based on the estimation of the second derivatives of the energy potential with respect to crack length using the technique of virtual domain perturbation $\boldsymbol{\theta}$~\cite{Dubois1998}. 
Numerically it uses an integral over a surface, which is more accurate than the contour integral used in the J-integral~\cite{rice1968plane}. 
Denoting the tangential vector to the crack tip as $\bm{t}$, we can use the virtual domain perturbation:
\begin{equation}
 \boldsymbol{\theta} = f(r) \bm{t}
 ,
\end{equation}
where
\begin{equation}
f(r) = \begin{cases}
1 & \; \text{for } r < r_i \\
\dfrac{r-r_o}{r_i-r_o} & \; \text{for } r_i < r < r_o \\
0 & \; \text{for } r_o < r
\end{cases}
\end{equation}
and $r$ is the distance from the crack tip, and $r_i$ and $r_o$ are set as $r_i=4\ell$ and $r_o=2.5r_i$ in this study (Figure~\ref{fig:Gtheta_schematic}).
With $\boldsymbol{\theta}$, the energy release rate is computed as
\begin{equation}
\label{eq:gtheta}
G_\theta = \int_\Omega \boldsymbol{\sigma} : ( \nabla \mathbf{u} \nabla \boldsymbol{\theta} )
- \frac{1}{2} (\boldsymbol{\sigma}: \nabla \mathbf{u} ) \nabla \cdot \boldsymbol{\theta} \, \mathd x
.
\end{equation}

\begin{figure}[htbp]
\centering %
\subfigure[ $G_{\theta}$]{\includegraphics[width=0.35\linewidth]{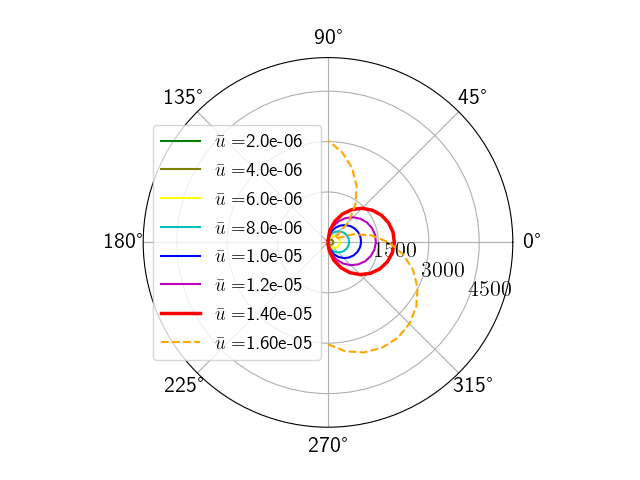}}
\quad
\subfigure[$G_{\theta}^t$]{\includegraphics[width=0.35\linewidth]{Figs/G_theta_Gcmax_polar_data_shear_Masonry_tensile_0.png}}
\caption{Comparison (a) energy release rate $G_{\theta}$ with (b) tensile part of energy release rate for \texttt{no-tension},  $\alpha=0$ under shear loading. In terms of directionality, tensile part of energy release rate accurately depicts the fracture propagation direction before and after fracture starts to propagate.}
\label{fig:G_theta_sigma_plus}
\end{figure}

However, the energy release rates computed by~\eqref{eq:gtheta} do not seem to predict the fracture propagation direction properly (Figure~\ref{fig:G_theta_sigma_plus} (a)). 
For this reason, considering the tension-compression asymmetry, we decompose $\mathbf{u}^t$ as
$$
\nabla \mathbf{u} = \nabla \mathbf{u}^t + \nabla \mathbf{u}^c
.
$$
Then we define the tensile part of the energy release rate as
\begin{equation}
\label{eq:gtheta_tens}
G_\theta^t = \int_\Omega \boldsymbol{\sigma}^t : ( \nabla \mathbf{u}^t \nabla \boldsymbol{\theta} )
- \frac{1}{2} (\boldsymbol{\sigma}^t: \nabla \mathbf{u}^t ) \nabla \cdot \boldsymbol{\theta} \, \mathd x
.
\end{equation}

While $\boldsymbol{\sigma}^t$ can be straightforwardly obtained from $\mathbb{C}\bm{\varepsilon}^t$, the tensile part of the displacement $\mathbf{u}^t$ is not defined in any of the strain energy split models.
Here, we assume that the compressive part of the displacement $\nabla \mathbf{u}^c$ is symmetric i.e., $(\nabla \mathbf{u}^c) = (\nabla \mathbf{u}^c)^\mathrm{T}$.
Then we have 
$$
\boldsymbol{\varepsilon}^c = \nabla \mathbf{u}^c
.
$$
and
\begin{equation}
\nabla \mathbf{u}^t 
= \nabla \mathbf{u} - \boldsymbol{\varepsilon}^c
.
\end{equation}
Figure~\ref{fig:G_theta_sigma_plus} illustrates the effects of the modified computation using the tensile part of the displacement, $G_\theta^t$ .
Without the modifications (Figure~\ref{fig:G_theta_sigma_plus} (a)), the maximum energy release rates point to an incorrect propagation direction (straight fracture propagation) prior to the propagation.
With our proposed modifications (Figure~\ref{fig:G_theta_sigma_plus} (b)), the direction predicted by the energy release rates coincide with the actual propagation direction.

    \bibliographystyle{apalike}
    \bibliography{main}

\begin{thebibliography}{}

\bibitem[Alessi et~al., 2018]{Alessi2017}
Alessi, R., Marigo, J.-J., Maurini, C., and Vidoli, S. (2018).
\newblock Coupling damage and plasticity for a phase-field regularisation of
  brittle, cohesive and ductile fracture: One--dimensional examples.
\newblock {\em International Journal of Mechanical Sciences}, 149:559--576.

\bibitem[Alessi et~al., 2014]{alessi2014gradient}
Alessi, R., Marigo, J.-J., and Vidoli, S. (2014).
\newblock Gradient damage models coupled with plasticity and nucleation of
  cohesive cracks.
\newblock {\em Archive for Rational Mechanics and Analysis}, 214(2):575--615.

\bibitem[Amadei, 1996]{Amadei1996}
Amadei, B. (1996).
\newblock Importance of anisotropy when estimating and measuring in situ
  stresses in rock.
\newblock {\em International Journal of Rock Mechanics and Mining Sciences \&
  Geomechanics Abstracts}, 33(3):293--325.

\bibitem[Ambati et~al., 2015]{Ambati2015}
Ambati, M., Gerasimov, T., and {De Lorenzis}, L. (2015).
\newblock {Phase-field modeling of ductile fracture}.
\newblock {\em Computational Mechanics}, 55(5):1017--1040.

\bibitem[Amor et~al., 2009]{Amor2009}
Amor, H., Marigo, J.-J., and Maurini, C. (2009).
\newblock Regularized formulation of the variational brittle fracture with
  unilateral contact: Numerical experiments.
\newblock {\em Journal of the Mechanics and Physics of Solids},
  57(8):1209--1229.

\bibitem[Baldelli et~al., 2014]{baldelli2014variational}
Baldelli, A.~L., Babadjian, J.-F., Bourdin, B., Henao, D., and Maurini, C.
  (2014).
\newblock A variational model for fracture and debonding of thin films under
  in-plane loadings.
\newblock {\em Journal of the Mechanics and Physics of Solids}, 70:320--348.

\bibitem[Bilgen and Weinberg, 2021]{bilgen2021phase}
Bilgen, C. and Weinberg, K. (2021).
\newblock A phase-field approach to pneumatic fracture with anisotropic crack
  resistance.
\newblock {\em International Journal of Fracture}, 232(2):135--151.

\bibitem[Bilke et~al., 2019]{Bilke2019337}
Bilke, L., Flemisch, B., Kalbacher, T., Kolditz, O., Helmig, R., and Nagel, T.
  (2019).
\newblock Development of open-source porous media simulators: Principles and
  experiences.
\newblock {\em Transport in Porous Media}, 130(1):337--361.

\bibitem[Bleyer and Alessi, 2018]{bleyer2018phase}
Bleyer, J. and Alessi, R. (2018).
\newblock Phase-field modeling of anisotropic brittle fracture including
  several damage mechanisms.
\newblock {\em Computer Methods in Applied Mechanics and Engineering},
  336:213--236.

\bibitem[Borden et~al., 2012]{Borden2012}
Borden, M., Verhoosel, C., Scott, M., Hughes, T., and Landis, C. (2012).
\newblock {A phase-field description of dynamic brittle fracture}.
\newblock {\em Computer Methods in Applied Mechanics and Engineering},
  217-220:77--95.

\bibitem[Bourdin et~al., 2012]{Bourdin2012}
Bourdin, B., Chukwudozie, C., and Yoshioka, K. (2012).
\newblock A variational approach to the numerical simulation of hydraulic
  fracturing.
\newblock In {\em SPE ATCE 2012}.

\bibitem[Bourdin et~al., 2000]{Bourdin2000}
Bourdin, B., Francfort, G., and Marigo, J.-J. (2000).
\newblock Numerical experiments in revisited brittle fracture.
\newblock {\em J. Mech. and Phys. of Solids}, 48(4):797--826.

\bibitem[Bourdin and Francfort, 2019]{Bourdin2019}
Bourdin, B. and Francfort, G.~A. (2019).
\newblock Past and present of variational fracture.
\newblock {\em SIAM News}, 52(9).

\bibitem[Bourdin et~al., 2011]{Bourdin2011}
Bourdin, B., Larsen, C., and Richardson, C. (2011).
\newblock {A time-discrete model for dynamic fracture based on crack
  regularization}.
\newblock {\em International Journal of Fracture}, 168(2):133--143.

\bibitem[Bryant and Sun, 2018]{bryant2018mixed}
Bryant, E.~C. and Sun, W. (2018).
\newblock A mixed-mode phase field fracture model in anisotropic rocks with
  consistent kinematics.
\newblock {\em Computer Methods in Applied Mechanics and Engineering},
  342:561--584.

\bibitem[Carcione, 1995]{carcione1995constitutive}
Carcione, J.~M. (1995).
\newblock Constitutive model and wave equations for linear, viscoelastic,
  anisotropic media.
\newblock {\em Geophysics}, 60(2):537--548.

\bibitem[Chambolle et~al., 2009]{Chambolle2009kink}
Chambolle, A., Francfort, G.~A., and Marigo, J.-J. (2009).
\newblock When and how do cracks propagate?
\newblock {\em Journal of the Mechanics and Physics of Solids},
  57(9):1614--1622.

\bibitem[De~Lorenzis and Maurini, 2021]{de2021nucleation}
De~Lorenzis, L. and Maurini, C. (2021).
\newblock Nucleation under multi-axial loading in variational phase-field
  models of brittle fracture.
\newblock {\em International Journal of Fracture}, pages 1--21.

\bibitem[Dubois et~al., 1998]{Dubois1998}
Dubois, F., Chazal, C., and Petit, C. (1998).
\newblock A finite element analysis of creep-crack growth in viscoelastic
  media.
\newblock {\em Mechanics Time-Dependent Materials}, 2(3):269--286.

\bibitem[Fei and Choo, 2021]{fei2021double}
Fei, F. and Choo, J. (2021).
\newblock Double-phase-field formulation for mixed-mode fracture in rocks.
\newblock {\em Computer Methods in Applied Mechanics and Engineering},
  376:113655.

\bibitem[Francfort and Marigo, 1998]{Francfort1998}
Francfort, G.~A. and Marigo, J.-J. (1998).
\newblock Revisiting brittle fracture as an energy minimization problem.
\newblock {\em Journal of the Mechanics and Physics of Solids},
  46(8):1319--1342.

\bibitem[Freddi and Royer-Carfagni, 2010]{Freddi2010}
Freddi, F. and Royer-Carfagni, G. (2010).
\newblock {Regularized variational theories of fracture: A unified approach}.
\newblock {\em Journal of the Mechanics and Physics of Solids},
  58(8):1154--1174.

\bibitem[Gerasimov and De~Lorenzis, 2022]{gerasimov2022second}
Gerasimov, T. and De~Lorenzis, L. (2022).
\newblock Second-order phase-field formulations for anisotropic brittle
  fracture.
\newblock {\em Computer Methods in Applied Mechanics and Engineering},
  389:114403.

\bibitem[Hakim and Karma, 2009]{Hakim2009}
Hakim, V. and Karma, A. (2009).
\newblock {Laws of crack motion and phase-field models of fracture}.
\newblock {\em Journal of the Mechanics and Physics of Solids}, 57(2):342--368.

\bibitem[He and Shao, 2019]{He-2019}
He, Q.-C. and Shao, Q. (2019).
\newblock {Closed-Form Coordinate-Free Decompositions of the Two-Dimensional
  Strain and Stress for Modeling Tension–Compression Dissymmetry}.
\newblock {\em Journal of Applied Mechanics}, 86(3).

\bibitem[Heider and Markert, 2017]{Heider2017}
Heider, Y. and Markert, B. (2017).
\newblock {A phase-field modeling approach of hydraulic fracture in saturated
  porous media}.
\newblock {\em Mechanics Research Communications}, 80:38--46.

\bibitem[Heng et~al., 2015]{heng2015experimental}
Heng, S., Guo, Y., Yang, C., Daemen, J.~J., and Li, Z. (2015).
\newblock Experimental and theoretical study of the anisotropic properties of
  shale.
\newblock {\em International Journal of Rock Mechanics and Mining Sciences},
  74:58--68.

\bibitem[Kuhn et~al., 2016]{Kuhn2016}
Kuhn, C., Noll, T., and M{\"{u}}ller, R. (2016).
\newblock {On phase field modeling of ductile fracture}.
\newblock {\em GAMM Mitteilungen}, 39(1):35--54.

\bibitem[Li and Maurini, 2019]{Li2019aniso}
Li, B. and Maurini, C. (2019).
\newblock Crack kinking in a variational phase-field model of brittle fracture
  with strongly anisotropic surface energy.
\newblock {\em Journal of the Mechanics and Physics of Solids}, 125:502--522.

\bibitem[Li et~al., 2015]{li2015phase}
Li, B., Peco, C., Mill{\'a}n, D., Arias, I., and Arroyo, M. (2015).
\newblock Phase-field modeling and simulation of fracture in brittle materials
  with strongly anisotropic surface energy.
\newblock {\em International Journal for Numerical Methods in Engineering},
  102(3-4):711--727.

\bibitem[Li et~al., 2021]{li2021research}
Li, C., Yang, D., Xie, H., Ren, L., and Wang, J. (2021).
\newblock Research on the anisotropic fracture behavior and the corresponding
  fracture surface roughness of shale.
\newblock {\em Engineering Fracture Mechanics}, 255:107963.

\bibitem[Li et~al., 2016]{Li2016grad}
Li, T., Marigo, J.-J., Guilbaud, D., and Potapov, S. (2016).
\newblock Gradient damage modeling of brittle fracture in an explicit dynamics
  context.
\newblock {\em International Journal for Numerical Methods in Engineering},
  108(11):1381--1405.

\bibitem[Luo et~al., 2022]{luo2022phase}
Luo, Z., Chen, L., Wang, N., and Li, B. (2022).
\newblock A phase-field fracture model for brittle anisotropic materials.
\newblock {\em Computational Mechanics}, pages 1--13.

\bibitem[Marigo et~al., 2016]{Marigo2016}
Marigo, J.-J., Maurini, C., and Pham, K. (2016).
\newblock An overview of the modelling of fracture by gradient damage models.
\newblock {\em Meccanica}, 51(12):3107--3128.

\bibitem[Mart{\'\i}nez-Pa{\~n}eda et~al., 2018]{martinez2018phase}
Mart{\'\i}nez-Pa{\~n}eda, E., Golahmar, A., and Niordson, C.~F. (2018).
\newblock A phase field formulation for hydrogen assisted cracking.
\newblock {\em Computer Methods in Applied Mechanics and Engineering},
  342:742--761.

\bibitem[Miehe et~al., 2016]{miehe2016ductile}
Miehe, C., Aldakheel, F., and Raina, A. (2016).
\newblock Phase field modeling of ductile fracture at finite strains: A
  variational gradient-extended plasticity-damage theory.
\newblock {\em International Journal of Plasticity}, 84:1--32.

\bibitem[Miehe et~al., 2010a]{Miehe2010}
Miehe, C., Hofacker, M., and Welschinger, F. (2010a).
\newblock A phase field model for rate-independent crack propagation: {R}obust
  algorithmic implementation based on operator splits.
\newblock {\em Computer Methods in Applied Mechanics and Engineering},
  199(45):2765--2778.

\bibitem[Miehe et~al., 2010b]{Miehe2010_variational}
Miehe, C., Welschinger, F., and Hofacker, M. (2010b).
\newblock {Thermodynamically consistent phase-field models of fracture:
  variational principles and multi-field FE implementations}.
\newblock {\em International Journal for Numerical Methods in Engineering},
  83(February):1273--1311.

\bibitem[Moerman et~al., 2016]{Moerman2016}
Moerman, K.~M., Simms, C.~K., and Nagel, T. (2016).
\newblock {Control of tension–compression asymmetry in Ogden hyperelasticity
  with application to soft tissue modelling}.
\newblock {\em Journal of the Mechanical Behavior of Biomedical Materials},
  56:218--228.

\bibitem[Nguyen et~al., 2020]{Nguyen2020}
Nguyen, T.-T., Yvonnet, J., Waldmann, D., and He, Q.-C. (2020).
\newblock Implementation of a new strain split to model unilateral contact
  within the phase field method.
\newblock {\em International Journal for Numerical Methods in Engineering},
  121(21):4717--4733.

\bibitem[Nguyen and Wu, 2018]{Nguyen2018}
Nguyen, V.~P. and Wu, J.-Y. (2018).
\newblock Modeling dynamic fracture of solids with a phase-field regularized
  cohesive zone model.
\newblock {\em Computer Methods in Applied Mechanics and Engineering},
  340:1000--1022.

\bibitem[Nixon et~al., 2010]{nixon2010anisotropic}
Nixon, M.~E., Cazacu, O., and Lebensohn, R.~A. (2010).
\newblock Anisotropic response of high-purity $\alpha$-titanium: Experimental
  characterization and constitutive modeling.
\newblock {\em International Journal of Plasticity}, 26(4):516--532.

\bibitem[Noii et~al., 2020]{noii2020adaptive}
Noii, N., Aldakheel, F., Wick, T., and Wriggers, P. (2020).
\newblock An adaptive global--local approach for phase-field modeling of
  anisotropic brittle fracture.
\newblock {\em Computer Methods in Applied Mechanics and Engineering},
  361:112744.

\bibitem[Parisio and Laloui, 2018]{parisio2018formulation}
Parisio, F. and Laloui, L. (2018).
\newblock On the formulation of anisotropic--polyaxial failure criteria: a
  comparative study.
\newblock {\em Rock Mechanics and Rock Engineering}, 51(2):479--489.

\bibitem[Pham et~al., 2011]{Pham2011Gradient}
Pham, K., Amor, H., Marigo, J.-J., and Maurini, C. (2011).
\newblock Gradient damage models and their use to approximate brittle fracture.
\newblock {\em International Journal of Damage Mechanics}, 20(4):618--652.

\bibitem[Rezaei et~al., 2022]{rezaei2022anisotropic}
Rezaei, S., Harandi, A., Brepols, T., and Reese, S. (2022).
\newblock An anisotropic cohesive fracture model: advantages and limitations of
  length-scale insensitive phase-field damage models.
\newblock {\em Engineering Fracture Mechanics}, page 108177.

\bibitem[Rice and Rosengren, 1968]{rice1968plane}
Rice, J.~R. and Rosengren, G. (1968).
\newblock Plane strain deformation near a crack tip in a power-law hardening
  material.
\newblock {\em Journal of the Mechanics and Physics of Solids}, 16(1):1--12.

\bibitem[Schuler et~al., 2020]{schuler2020chemo}
Schuler, L., Ilgen, A.~G., and Newell, P. (2020).
\newblock Chemo-mechanical phase-field modeling of dissolution-assisted
  fracture.
\newblock {\em Computer Methods in Applied Mechanics and Engineering},
  362:112838.

\bibitem[Steinke and Kaliske, 2019]{Steinke2018}
Steinke, C. and Kaliske, M. (2019).
\newblock A phase-field crack model based on directional stress decomposition.
\newblock {\em Computational Mechanics}, 63(5):1019--1046.

\bibitem[Teichtmeister et~al., 2017]{teichtmeister2017phase}
Teichtmeister, S., Kienle, D., Aldakheel, F., and Keip, M.-A. (2017).
\newblock Phase field modeling of fracture in anisotropic brittle solids.
\newblock {\em International Journal of Non-Linear Mechanics}, 97:1--21.

\bibitem[Ulloa et~al., 2022]{ulloa2022micromechanics}
Ulloa, J., Wambacq, J., Alessi, R., Samaniego, E., Degrande, G., and
  Fran{\c{c}}ois, S. (2022).
\newblock A micromechanics-based variational phase-field model for fracture in
  geomaterials with brittle-tensile and compressive-ductile behavior.
\newblock {\em Journal of the Mechanics and Physics of Solids}, 159:104684.

\bibitem[van Dijk et~al., 2020]{DIJK2020}
van Dijk, N., Espadas-Escalante, J., and Isaksson, P. (2020).
\newblock Strain energy density decompositions in phase-field fracture theories
  for orthotropy and anisotropy.
\newblock {\em International Journal of Solids and Structures},
  196-197:140--153.

\bibitem[Wang et~al., 2019]{wang2019anisotropic}
Wang, Y., Tan, W., Liu, D., Hou, Z., and Li, C. (2019).
\newblock On anisotropic fracture evolution and energy mechanism during marble
  failure under uniaxial deformation.
\newblock {\em Rock Mechanics and Rock Engineering}, 52(10):3567--3583.

\bibitem[Wheeler et~al., 2014]{Wheeler2014}
Wheeler, M., Wick, T., and Wollner, W. (2014).
\newblock {An augmented-Lagrangian method for the phase-field approach for
  pressurized fractures}.
\newblock {\em Computer Methods in Applied Mechanics and Engineering},
  271:69--85.

\bibitem[Wilson and Landis, 2016]{Wilson2016}
Wilson, Z. and Landis, C. (2016).
\newblock {Phase-field modeling of hydraulic fracture}.
\newblock {\em Journal of the Mechanics and Physics of Solids}, 96:264--290.

\bibitem[Wu and Nguyen, 2018]{wu2018length}
Wu, J.-Y. and Nguyen, V.~P. (2018).
\newblock A length scale insensitive phase-field damage model for brittle
  fracture.
\newblock {\em Journal of the Mechanics and Physics of Solids}, 119:20--42.

\bibitem[Yin and Kaliske, 2020]{Yin2020}
Yin, B. and Kaliske, M. (2020).
\newblock A ductile phase-field model based on degrading the fracture
  toughness: {T}heory and implementation at small strain.
\newblock {\em Computer Methods in Applied Mechanics and Engineering},
  366:113068.

\bibitem[Yoshioka et~al., 2021]{Yoshioka2021int}
Yoshioka, K., Mollaali, M., and Kolditz, O. (2021).
\newblock Variational phase-field fracture modeling with interfaces.
\newblock {\em Computer Methods in Applied Mechanics and Engineering},
  384:113951.

\bibitem[Yoshioka et~al., 2019]{Yoshioka2019}
Yoshioka, K., Parisio, F., Naumov, D., Lu, R., Kolditz, O., and Nagel, T.
  (2019).
\newblock {Comparative verification of discrete and smeared numerical
  approaches for the simulation of hydraulic fracturing}.
\newblock {\em GEM - International Journal on Geomathematics}, 10(1).

\bibitem[You et~al., 2021]{you2021brittle}
You, T., Waisman, H., and Zhu, Q.-Z. (2021).
\newblock Brittle-ductile failure transition in geomaterials modeled by a
  modified phase-field method with a varying damage-driving energy coefficient.
\newblock {\em International Journal of Plasticity}, 136:102836.

\bibitem[Zhou et~al., 2018]{zhou2018phase}
Zhou, S., Zhuang, X., and Rabczuk, T. (2018).
\newblock A phase-field modeling approach of fracture propagation in
  poroelastic media.
\newblock {\em Engineering Geology}, 240:189--203.

\end{thebibliography}
\end{document}